\documentclass[10pt]{amsart}

\let\oldmarginpar\marginpar
\renewcommand\marginpar[1]{\-\oldmarginpar[\raggedleft\footnotesize #1]
{\raggedright\footnotesize #1}}

\DeclareMathSizes{6}{18}{12}{8}

\usepackage{mathptmx}
\usepackage{lineno}
\usepackage{amsmath}
\usepackage{amscd}
\usepackage{amssymb}
\usepackage{amsthm}
\usepackage{xspace}
\usepackage[all,tips]{xy}
\usepackage[dvips]{graphicx}
\usepackage{verbatim}
\usepackage{syntonly}
\usepackage{hyperref}
\usepackage{color}

\author{Benjamin Linowitz}
\address{Department of Mathematics\\University of Michigan\\Ann Arbor, MI 48109}
\email{linowitz@umich.edu}

\author{D. B. McReynolds}
\address{Department of Mathematics\\Purdue University\\West Lafayette, IN 47907}
\email{dmcreyno@math.purdue.edu}

\author{Nicholas Miller}
\address{Department of Mathematics\\Purdue University\\West Lafayette, IN 47907}
\email{mille965@math.purdue.edu}

\theoremstyle{plain}
\newtheorem{thm}{Theorem}[section]
\newtheorem{cor}[thm]{Corollary}
\newtheorem{prop}[thm]{Proposition}
\newtheorem{lemma}[thm]{Lemma}

\DeclareMathOperator{\Aut}{Aut} 
 
\DeclareMathOperator{\SL}{SL} \DeclareMathOperator{\PSL}{PSL}

 \DeclareMathOperator{\Ad}{Ad}
 \DeclareMathOperator{\Mat}{Mat}

\DeclareMathOperator{\FI}{FI}

\DeclareMathOperator{\D}{D}\DeclareMathOperator{\Ram}{Ram}
\DeclareMathOperator{\Br}{Br}
\DeclareMathOperator{\Gal}{Gal}

\DeclareMathOperator{\Res}{Res}
\DeclareMathOperator{\Ord}{Ord}
\DeclareMathOperator{\Inv}{Inv}
\DeclareMathOperator{\Spec}{Spec}\DeclareMathOperator{\ord}{ord}
\DeclareMathOperator{\V}{V}\DeclareMathOperator{\Para}{Para}
\DeclareMathOperator{\Build}{Build}

\DeclareMathOperator{\Aa}{A}

\newcommand{\vp}{\varphi}

\newcommand{\iny}{\infty}

\newcommand{\abs}[1]{\left\vert#1\right\vert}

\newcommand{\ceil}[1]{\left\lceil #1\right\rceil}
\newcommand{\set}[1]{\left\{#1\right\}}

\newcommand{\pr}[1]{\left( #1 \right) }
\newcommand{\su}{\subset}

\newcommand{\lra}{\longrightarrow}

\newcommand{\B}[1]{\ensuremath{\mathbf{#1}}}

\newcommand{\Cal}[1]{\ensuremath{\mathcal{#1}}}

\newcommand{\N}{\ensuremath{\mathbf{N}}}
\newcommand{\Q}{\ensuremath{\mathbf{Q}}}
\newcommand{\R}{\ensuremath{\mathbf{R}}}
\newcommand{\Z}{\ensuremath{\mathbf{Z}}}
\newcommand{\C}{\ensuremath{\mathbf{C}}}
\newcommand{\spec}[1]{\Spec(#1)}

\title{Locally Equivalent Correspondences}

\begin{document}


\maketitle

\begin{abstract}
\noindent Given a pair of number fields with isomorphic rings of adeles, we construct bijections between objects associated to the pair. For instance we construct an isomorphism of Brauer groups that commutes with restriction. We additionally construct bijections between central simple algebras, maximal orders, various Galois cohomology sets, and commensurability classes of arithmetic lattices in simple, inner algebraic groups. We show that under certain conditions, lattices corresponding to one another under our bijections have the same covolume and pro-congruence completion. We also make effective a finiteness result of Prasad and Rapinchuk.
\end{abstract}

\section{Introduction}

Given a number field $K$, we denote the Brauer group of $K$ by $\Br(K)$. For any subfield $F \subset K$, we have a homomorphism $\Res_{K/F}\colon \Br(F) \to \Br(K)$ given by $\Res_{K/F}([B]) = [B \otimes_F K]$. For a pair of number fields $K,K'$, a \textbf{natural isomorphism} between $\Br(K),\Br(K')$ is an isomorphism $\Phi_{\Br}\colon \Br(K) \to \Br(K')$ such that for any $F \su K \cap K'$ and any $L$ with $KK' \su L$, the diagram
\begin{equation}\label{BrauerNatural} 
\xymatrixrowsep{.2in} \xymatrixcolsep{.2in} \xymatrix{&\Br(L)&\\\Br(K)\ar^{\Res}[ru]\ar@{<->}[rr]^{\Phi_{\Br}}&&\Br(K')\ar[lu]_{\Res}\\&\Br(F)\ar[ru]_\Res \ar[lu]^\Res&}
\end{equation}
commutes; note that $\Br(K) \cong \Br(L)$ as abstract groups provided $K,L$ have the same number of real places. The fiber or pullback of a class $[A] \in \Br(K)$ under the map $\Res_{K/F}$ gives a family of $F$--subalgebras of $A$. In \cite{McReid} (see also \cite{LMPT}), it was shown that these fibers determine the algebra in certain situations. As in \cite{McReynolds} however, there are situations when these fibers fail to determine the algebra. In particular, when a natural isomorphism between Brauer groups exists, the pair $[A],\Phi_{\Br}([A])$ provide examples for any $[A] \in \Br(K)$.

To construct natural isomorphisms we will make use of what we call locally equivalent number fields. For a number field $K$, denote by $\V^K$ the set of places of $K$ and by $\B{A}_K$ the ring of $K$--adeles. We say that $K$ and $K'$ are \textbf{locally equivalent} if there exists a bijection $\Phi_{\V}\colon \V^K \to \V^{K'}$ between places such that $K_v \cong K'_{\Phi_{\V}(v)}$ for all $v\in \V^K$. By work of Iwasawa \cite{Iwasawa}, this condition is equivalent to the condition that the two fields have isomorphic rings of adeles. We will refer to the pair $K,K'$ as a \textbf{locally equivalent pair} when $K,K'$ are locally equivalent number fields.

\begin{thm}\label{BrauerIsomorphism}
For any locally equivalent pair $K,K'$, there is a natural isomorphism $\Phi_{\Br} \colon \Br(K)\to \Br(K')$.
\end{thm}

We note that it is known that arbitrarily large families of pairwise locally equivalent, non-isomorphic number fields exist (see \cite{Komatsu3}). Locally equivalent fields or variants have recently been employed by Aka \cite{Aka} and D. Prasad \cite{DPra}. Aka used them to produce examples of incommensurable arithmetic lattices with the same profinite completions. D. Prasad used a refinement of arithmetic equivalence to produce Riemann surfaces with the same Jacobians viewed as complex abelian varieties. 

A natural isomorphism between Brauer groups induces a bijection between algebras in each Morita class. In order to pass to finer structures like lattices, we require a refinement of Theorem \ref{BrauerIsomorphism}. For a central simple algebra $A$ over $K$, we denote by $\Ord(A,K)$ the set of $\Cal{O}_K$--orders of $A$ having full rank. Our next result exhibits a close relationship between the set of orders of $A$ and of $\Phi_{\Br}(A)$; we refer the reader to Subsection \ref{OrderSubSec} for the definition of the level ideal. 

\begin{thm}\label{orderbij} 
If $K,K'$ are a locally equivalent pair, then there exists a bijection 
\[ \Phi_{\Ord}\colon \Ord(A,K) \lra \Ord(\Phi_{\Br}(A),K') \] 
with the property that if $\Cal{R}\in \Ord(A,K)$ has level ideal $\Cal{L}_{\Cal{R}}$ then $\Phi_{\Ord}(\Cal{R})$ has level ideal $\Phi_{\V}(\Cal{L}_{\Cal{R}})$.
\end{thm}

We will employ Theorem \ref{orderbij} to establish a bijection between maximal arithmetic lattices arising from central simple algebras defined over locally equivalent fields. This bijection will be denoted by $\Phi_{\mathrm{lattice}}$. Our next result shows that covolume and pro-congruence topology are preserved under our bijection. We note that if $\Lambda$ is a lattice for which we can apply our bijection and $G$ is the associated semisimple Lie group, then $\Phi_{\mathrm{lattice}}(\Lambda)$ is also a lattice in the semisimple Lie group $G$.

\begin{cor}\label{PropertyBiject}
Under the bijection $\Phi_{\mathrm{lattice}}$, the lattices $\Lambda,\Phi_{\mathrm{lattice}}(\Lambda)$ have the same pro-congruence completion. If $\Lambda$ is derived from an order, then $\Lambda,\Phi_{\mathrm{lattice}}(\Lambda)$ have the same covolume.
\end{cor}

The bijection $\Phi_{\Br}$ is constructed using the methods of \cite{McReynolds}. The integral refinements to orders and maximal lattices are proven using a similar principle and make crucial use of the local-to-global correspondence for lattices. The relationship between pro-congruence topologies follows immediately from the construction of the bijection, as does the preservation of volume. The main tools required here are technical but standard and include Bruhat--Tits theory \cite{Tits}, the work of Borel \cite{borel} and Borel--Prasad \cite{BP} on the classification of maximal arithmetic lattices, and Prasad's volume formula \cite{Prasad}.

\subsection{Partial converses}

The converse of Theorem \ref{BrauerIsomorphism} is false in general as there exist arithmetically equivalent number fields with naturally isomorphic Brauer groyps which are not locally equivalent (see \cite{McReynolds}). We say that $K,K'$ are \textbf{locally GCD equivalent} if, for every rational prime $p$ which is unramified in both $K$ and $K'$, we have
\begin{equation}\label{blah} 
\gcd\set{ [K_v:\Q_p]~:~v\in V^{K},~v\mid p} = \gcd\set{[K'_{v'}:\Q_p]~:~v'\in V^{K'},~ v'\mid p}.
\end{equation}

\begin{thm}\label{GCDEQUIV}
If $K,K'$ are number fields for which there is a natural isomorphism between $\Br(K)$ and $\Br(K')$, then $K$ and $K'$ are locally GCD equivalent. 
\end{thm}

It is known that number fields which are either arithmetically equivalent or locally equivalent must have the same degree, discriminant, and Galois closure (\cite[Thm 1]{Perlis}). We do not know if the same is true for locally GCD equivalent fields. A straightforward application of Theorem \ref{GCDEQUIV} is the following rigidity result which generalizes \cite[Thm 1.1]{McReid} and  \cite[Thm 1.2]{LMPT}.

\begin{cor}\label{RIGIDITY}
If $K,K'$ are finite Galois extensions of $\Q$ and there is a natural isomorphism between $\Br(K)$ and $\Br(K')$, then $K \cong K'$.
\end{cor}

\subsection{Effective result of Prasad and Rapinchuk}

Our proof of Corollary \ref{PropertyBiject} will require us to examine a special case of Prasad's volume formula pertaining to central simple algebras (see Proposition \ref{P:volumeformula} below). Our next result is a further application of this formula and concerns a problem arising in spectral geometry.

Given a semisimple Lie group $G$ and a maximal compact subgroup $K$, we have an associated symmetric space $\mathcal{X}_G=G/K$. Any lattice $\Gamma$ in $G$ will give rise to a locally symmetric orbifold $M = \Gamma \setminus \mathcal{X}_G$. The \textbf{geodesic length spectrum} $\Cal{L}(M)$ of $M$ is the set of lengths of closed geodesics counted with multiplicity while the \textbf{geodesic length set} $\textrm{L}(M)$ is the set of lengths without multiplicity. We say two manifolds $M,N$ are \textbf{geodesic length isospectral} if $\Cal{L}(M)=\Cal{L}(N)$. We say two manifolds $M,N$ are \textbf{length commensurable} if $\Q\textrm{L}(M) = \Q\textrm{L}(N)$. Reid \cite{Reid} proved that if $M,N$ are finite volume hyperbolic 2--manifolds such that $M$ is arithmetic and $M,N$ are length commensurable, then $M,N$ are commensurable. In particular, $N$ must also be arithmetic. The result was extended to hyperbolic 3--manifolds by Chinburg--Hamilton--Long--Reid \cite{CHLR}. However, even before \cite{CHLR}, it was known that length commensurability does not imply commensurability in general. Lubotzky--Samuel--Vishne \cite{LSV} produced examples of incommensurable arithmetic lattices in $\SL(n,\R),\SL(n,\C)$ for all $n>2$ that are length isospectral. Prasad--Rapinchuk \cite{PrasadRap} generalized these works addressing precisely when the above commensurability rigidity holds; the most general versions rely on Schanuel's conjecture in transcendental number theory. They proved that for a fixed manifold $M$ of the above type, there are only finitely many commensurability classes of manifolds that can be length commensurable to $M$. Our next result provides an explicit upper bound for the number of classes as a function of only the volume of the manifold $M$. The class of manifolds $M$ are those arising from groups of the form $\SL_1(D)$, where $D$ is a division algebra defined over a number field.

\begin{thm}\label{EffectivePR}
If $K$ is a number field, $D$ a $K$--division algebra of degree $d>1$, $\mathcal{R}$ a maximal order of $D$, $\Gamma= \Res_{K/\Q}(\mathcal{R}^1)$ the arithmetic lattice in the semisimple Lie group $G=\Res_{K/\Q}(D^1)(\R)$, $M = \Gamma\setminus \mathcal{X}_G$, and $V = \mathrm{Vol}(M)$, then the number of pairwise non-commensurable manifolds that are length commensurable with $M$ is bounded above by $1+10^{33}V$.
\end{thm}

We point out that for division algebras of sufficiently large degree (in fact $d>28$) this bound reduces to a more aesthetically pleasing $1+V$. This strengthening is an immediate consequence of the proof. It is well-known (see for instance \cite[Thm 10.1]{PrasadRap}) that $\Cal{L}(M)$ is directly related to the eigenvalue spectrum $\mathcal{E}(\Delta_M)$ of the Laplace--Beltrami operator $\Delta_M$ acting on $L^2(M)$. Specifically, the Laplace--Beltrami spectrum determines $\textrm{L}(M)$. From Theorem \ref{EffectivePR}, we obtain the following corollary.

\begin{cor}
If $K$ is a number field, $D$ a $K$--division algebra of degree $d>1$, $\mathcal{R}$ a maximal order of $D$, $\Gamma= \Res_{K/\Q}(\mathcal{R}^1)$ the arithmetic lattice in the semisimple Lie group $G=\Res_{K/\Q}(D^1)(\R)$, $M = \Gamma\setminus \mathcal{X}_G$, and $V = \mathrm{Vol}(M)$, then the number of pairwise non-commensurable manifolds that are isospectral with $M$ is bounded above by $1+10^{33}V$.
\end{cor} 

\subsection{Galois cohomology sets and maximal lattices}

Returning to locally equivalent correspondences, the above correspondences between Brauer groups and maximal arithmetic lattices have extensions to other classes of arithmetic lattices and manifolds. The role of the Brauer group is played by Galois cohomology sets that parameterize commensurability classes of arithmetic manifolds. In fact, the Brauer group parameterizes inner forms of type $\Aa_{n-1}$ via its subgroup of $n$--torsion $\Br_n(K)$. For locally equivalent fields, we show these natural bijections also hold for Galois cohomology sets associated to absolutely almost simple, inner forms of split algebraic group.

\begin{thm}\label{GC-Isomorphisms}
For any locally equivalent pair $K,K'$ and any absolutely almost simple, split group $G$, there is a natural bijection between $H^1(K,\overline{G})$ and $H^1(K',\overline{G})$.
\end{thm}

As a consequence, the above theorem additionally holds for any absolutely almost simple groups which are inner forms of a given split group.  This is because their Galois cohomology sets are in natural bijection with those of Theorem \ref{GC-Isomorphisms}. As with the special case of Brauer groups, we also have bijections between maximal lattices in the commensurability classes.

\begin{thm}\label{GC-MaximalBiject}
If $[\xi']\in H^1(K',\overline{G})$ is the equivalence class of cocycles corresponding to $[\xi]\in H^1(K,\overline{G})$ in the bijection of Theorem \ref{GC-Isomorphisms} and $_\xi G$ and $_{\xi'}G'$ are the corresponding isomorphism classes of twists of $G$ and $G'$, then there is a bijection between maximal arithmetic lattices of $_\xi G$ and $_{\xi'}G'$.
\end{thm}

Under the bijections given by Theorem \ref{GC-MaximalBiject}, the associated pair of maximal lattices have the same pro-congruence completion. When the lattices also have the congruence subgroup property, one obtains non-isomorphic lattices with isomorphic profinite completions. These examples are not new as they appeared in \cite{Aka}. The volume is not a purely local invariant and volumes of associated manifolds under our bijection do not always agree; see the example at the end of Section 3. However, under our bijection, there are many cases when it does (e.g. Corollary \ref{PropertyBiject}). 

\paragraph{\textbf{Acknowledgements.}}

The authors would like to thank Ted Chinburg, Britain Cox, Rachel Davis, Bart de Smit, Matt Emerton, Amir Mohammadi, and Alan Reid for useful conversations on the material of this paper. We would also like to thank the referee for helpful comments and suggestions for revisions of an earlier draft of this paper. The first author was partially supported by an NSF RTG grant DMS-1045119 and an NSF Mathematical Sciences Postdoctoral Fellowship. The second author was partially supported by the NSF grants DMS-1105710 and DMS-1408458.

\section{Local Equivalence and Brauer Groups}\label{BrauerIsoSection}

In this section, we work out the details of Theorem \ref{BrauerIsomorphism} and the various refinements given in the introduction. 

\subsection{Brauer groups}

We begin by recalling some basic properties of the Brauer groups associated to local and global fields. We refer the reader to \cite{Pierce, Reiner} for a more detailed treatment. For a number field $K$, the \textbf{Brauer group} $\Br(K)$ is the group of Morita equivalence classes of central simple algebras defined over $K$ with the tensor product as the group operation. Given an extension of fields $L/K$, there is a well-defined restriction homomorphism $\Res_{L/K}\colon \Br(K) \to \Br(L)$ given by $\Res_{L/K}([B]) = [B\otimes_K L]$. Every Morita equivalence class contains a unique division algebra \cite[p. 228]{Pierce}. Consequently, when $B$ is a division algebra we will often simply write $\Res_{L/K}(B)=B\otimes_K L$. The Hasse invariants associated to the places $v$ of $K$ are defined as follows. For a finite place $v$ of $K$, the homomorphism $\Inv_v\colon \Br(K_v)\to \bf{Q}/\bf{Z}$ which sends a Brauer class to its associated Hasse invariant is an isomorphism \cite[p.~338]{Pierce}. It follows that the order of a class $[B_v]$ in $\Br(K_v)$ with Hasse invariant $\frac{a_v}{m_v}$ is equal to $m_v$. Here $a_v,m_v$ are non-negative relatively prime integers with $a_v\leq m_v$. For a complex archimedean place $v$ of $K$, any central simple algebra over $K_v$ is isomorphic to $\Mat(n,\bf{C})$ for some positive integer $n$ and consequently the group $\Br(K_v)$ is trivial. We define the local Hasse invariant to be 0 in this case. If $v$ is a real place then $\Br(K_v)\cong\bf{Z}/2\bf{Z}$. The latter group is generated by the equivalence class $[\bf{H}]$ of Hamilton's quaternions. In this case we define $\Inv_v([\textbf{R}])=0$ and $\Inv_v([\textbf{H}])=\frac{1}{2}$. It is a consequence of class field theory that the sequence
\begin{equation}\label{ABHN}
\xymatrix{ 1 \ar[r] & \Br(K) \ar[r]^-{\iota_{\Br,K}}  & \bigoplus_{v \in V^K} \Br(K_v) \ar[r]^-{\sigma} & \Q/\Z \ar[r] & 1},
\end{equation}
where $\iota_{\Br,K}(B)=\set{B_v}_v$ and $\sigma(\set{B_v}_v)=\sum_v \Inv_v(B_v)$, is exact (see \cite[Ch.~32]{Reiner}). The finite set of places of $K$ where $\Inv_v([B]) \ne 0$ is denoted by $\Ram([B])$ and the subsets of non-archimedean and archimeadean places are denoted by $\Ram_f([B])$ and $\Ram_\iny([B])$, respectively. 

\subsection{Proof of Theorem \ref{BrauerIsomorphism}}\label{BrauerIsoSubSec}

Recall that we have a pair of locally equivalent fields $K,K'$ and seek to produce a natural isomorphism $\Phi_{\Br}\colon \Br(K) \to \Br(K')$. To be a natural isomorphism, we require that $\Phi_{\Br}$ be a group isomorphism and additionally that whenever $F \su K,K' \su L$, commutativity of the following diagram:
\[ \xymatrixrowsep{.2in} \xymatrixcolsep{.2in}\xymatrix{&\Br(L)&\\\Br(K)\ar^{\Res}[ru]\ar@{<->}[rr]^{\Phi_{\Br}}&&\Br(K')\ar[lu]_{\Res}\\&\Br(F)\ar[ru]_\Res \ar[lu]^\Res&} \] 
Of course one knows the if $K$ and $K'$ have the same number of real place then their Brauer groups are isomorphic, however the latter condition does not follow in general.
We briefly remark that our construction of such an isomorphism depends on the choice of $\Phi_V$.

\begin{proof}[Proof of Theorem \ref{BrauerIsomorphism}]
As every Morita equivalence class in $\Br(K)$ is represented by a unique division algebra, it suffices to define $\Phi_{\Br}$ on the level of division algebras. To that end, we will make extensive use of the bijection $\Phi_{\V}\colon \V^K \to \V^{K'}$ mentioned in the introduction. Recall that under this bijection we have $K_v \cong K'_{\Phi_V(v)}$ for all $v \in V$. We define $\Phi_{\Br}\colon \Br(K) \to \Br(K')$ by the Hasse invariant equations
\begin{equation}\label{B2}
\Inv_{v'}(\Phi_{\Br}(A)) = \Inv_{\Phi_{\V}^{-1}(v')}(A).
\end{equation}
It follows by (\ref{ABHN}) that there exists a unique division algebra $A' = \Phi_{\Br}(A)$ with these Hasse invariants. We additionally have an inverse process for producing a map $\Phi_{\Br}^{-1}\colon \Br(K') \to \Br(K)$ given by the Hasse invariant equations
\begin{equation}\label{B4}
\Inv_{v}(\Phi_{\Br}^{-1}(A')) = \Inv_{\Phi_{\V}(v)}(A').
\end{equation}
Now given $[A] \in \Br(K)$, we see that $\Inv_{v'}(\Phi_{\Br}(A)) = \Inv_{\Phi_{\V}^{-1}(v')}(A)$ and that
\[ \Inv_{v}(\Phi_{\Br}^{-1}(\Phi_{\Br}(A))) = \Inv_{\Phi_{\V}(v)}(\Phi_{\Br}(A)) = \Inv_{\Phi_{\V}^{-1}(\Phi_{\V}(v))}(A) = \Inv_v(A). \]
Thus $\Inv_v(\Phi_{\Br}^{-1}(\Phi_{\Br}(A))) = \Inv_v(A)$, and hence by (\ref{ABHN}), $\Phi_{\Br}$ is a bijection. To see that $\Phi_{\Br}$ is a group homomorphism, simply note that $\Inv_v(A_1 \otimes_K A_2) = \Inv_v(A_1)+\Inv_v(A_2)$. It remains to show that the isomorphism $\Phi_{\Br}$ is natural. We first check that $\Phi_{\Br} \circ \Res_{K/F} = \Res_{K'/F}$. For each $[B]$ in $\Br(F)$, we must show that $\Phi_{\Br}([B \otimes_F K]) = [B \otimes_F K']$. By (\ref{ABHN}), it suffices to check that for each $v' \in \V^{K'}$, we have
\begin{equation}\label{B1}
\Inv_{v'}(\Phi_{\Br}([B \otimes_F K])) = \Inv_{v'}([B \otimes_F K']).
\end{equation}
Via the Hasse invariant equations (\ref{B2}), (\ref{B4}), we have $\Inv_v([B \otimes_F K]) = \Inv_{v'}(\Phi_{\Br}([B \otimes_F K]))$. Additionally, we have the equality $\Inv_v([B \otimes_F K]) = \Inv_{v'}([B \otimes_F K'])$, which follows from the local equivalence of $K,K'$ and basic properties of Hasse invariants. In total, we see that (\ref{B1}) holds. The verification of $\Phi_{\Br}^{-1} \circ \Res_{K'/F} = \Res_{K/F}$ is identical. For the top triangle in (\ref{BrauerNatural}), we first check that $\Res_{L/K'} \circ \Phi_{\Br} = \Res_{L/K}$. As in the first part, it suffices by (\ref{ABHN}) to verify
\begin{equation}\label{B3}
\Inv_w([A \otimes_K L]) = \Inv_w([\Phi_{\Br}(A) \otimes_{K'} L]).
\end{equation}
As in the first case, (\ref{B3}) follows from from the local equivalence of $K,K'$ in combination with basic properties of Hasse invariants and equations (\ref{B2}), (\ref{B4}). The verification of $\Res_{L/K'} = \Res_{L/K} \circ \Phi_{\Br}^{-1}$ is identical.
\end{proof} 

\subsection{Proof of Theorem \ref{orderbij}}\label{OrderSubSec}

We denote the extension of $\Phi_{\V}$ to a bijection between fractional ideals by $\Phi_{\FI}$. Following \cite[p. 49]{Reiner}, if $X$ is a finitely generated $\mathcal{O}_K$--module, we define the \textbf{order ideal}, denoted $\ord(X)$, by the convention:
\begin{enumerate}
\item[(a)] If $X=0$, then $\ord(X) = \mathcal{O}_K$.
\item[(b)] If $X$ is not an $\mathcal{O}_K$--torsion module, $\ord(X)=0$.
\item[(c)] If $X$ is a nonzero $\mathcal{O}_K$--torsion module, it has an $\mathcal{O}_K$--composition series with factors $\{\mathcal{O}_K/p_i\}$ where $p_i$ ranges over some set of maximal ideals.  We then set $\ord(X)=\prod p_i$.
\end{enumerate}

For a central simple algebra $A$ over $K$, we denote by $\Ord(A,K)$ the set of $\Cal{O}_K$--orders of $A$. Given an order $\Cal{R}$ in $\Ord(A,K)$, we define the \textbf{level ideal} $\Cal{L}_{\Cal{R}}$ of $\Cal{R}$ to be the order ideal $\textrm{ord}(\mathcal{O}/\mathcal{R})$ of the $\Cal{O}_K$--module $\Cal{O}/\Cal{R}$, where $\Cal{O}$ is a maximal order of $A$ containing $\Cal{R}$. This definition is independent of the choice of $\Cal{O}$. We remark that as an immediate consequence of the definition, $\mathcal{O}$ is a maximal order in $A$ if and only if $\mathcal{L}_\mathcal{O}=\text{ord }(0)=\mathcal{O}_K$.

\begin{proof}[Proof of Theorem \ref{orderbij}]
Let $\Cal{O}$ be a maximal order of $A$ and $\Cal{O}'$ be a maximal order of $\Phi_{\Br}(A)$. By Theorem \ref{BrauerIsomorphism} we know that for every place $v\in \V^K$ there is an isomorphism $\Phi_v\colon  A\otimes_K K_v \to \Phi_{\Br}(A)\otimes_{K'} K'_{\Phi_{\V}(v)}$. Since all of the maximal orders of $\Phi_{\Br}(A)\otimes_{K'} K'_{\Phi_{\V}(v)}$ are conjugate, we may assume without loss of generality that $\Phi_v(\Cal{O}\otimes_{\Cal{O}_K}\Cal{O}_{K_v})=\Cal{O}'\otimes_{\Cal{O}_{K'}} \Cal{O}_{K'_{\Phi_{\V}(v)}}$. We now define the map $\Phi_{\Ord}$. For an arbitrary $\Cal{O}_K$--order $\Cal{R}$ of $A$, we define $\Phi_{\Ord}(\Cal{R})$ to be the unique $\Cal{O}_{K'}$--order of $\Phi_{\Br}(A)$ whose completions satisfy the equation
\begin{equation}\label{C1}
\Phi_{\Ord}(\Cal{R})\otimes_{\Cal{O}_{K'}} \Cal{O}_{K'_{\Phi_{\V}(v)}} = \Phi_v(\Cal{R}\otimes_{\Cal{O}_K} \Cal{O}_{K_v})
\end{equation}
for each place $v$ in $\V^K$. That such an order $\Phi_{\Ord}(\Cal{R})$ exists follows from the local-to-global correspondence for orders (see \cite[Thm 4.22]{Reiner} for instance). Indeed, it suffices to show that there is an $\Cal{O}_{K'}$--order $\Cal{R}'$ of $\Phi_{\Br}(A)$ such that $\Cal{R}'\otimes_{\Cal{O}_{K'}} \Cal{O}_{K'_{\Phi_{\V}({v})}}=\Phi_v(\Cal{R}\otimes_{\Cal{O}_K} \Cal{O}_{K_v})$ for all but finitely many places $v\in \V^K$. As $\Cal{O}\otimes_{\Cal{O}_K}\Cal{O}_{K_v}=\Cal{R}\otimes_{\Cal{O}_K}\Cal{O}_{K_v}$ for all but finitely places $v\in \V^K$, it is clear that $\Cal{O}'$ has the required property.  We now show that $\Phi_{\Ord}$ is surjective, as the injectivity of $\Phi_{\Ord}$ is clear. To that end, let $\Cal{R}'\in\Ord(\Phi_{\Br}(A),K')$ and $\{v_1,\dots,v_n \}$ be the set of places of $\V^{K'}$ for which $\Cal{R}'\otimes_{\Cal{O}_{K'}} \Cal{O}_{K'_{v_i}} \ne \Cal{O}'\otimes_{\Cal{O}_{K'}} \Cal{O}_{K'_{v_i}}$. Also, set $\Cal{R}$ to be the $\Cal{O}_K$--order of $A$ whose completions are equal to $\Cal{O}\otimes_{\Cal{O}_K} \Cal{O}_{K_v}$ if $\Phi_{\V}(v)\not\in \{ v_1,\dots,v_n \}$. Otherwise, we set the completion to have image in $\Phi_{\Br}(A)\otimes_{K'} K'_{v_i}$ equal to $\Cal{R}'\otimes_{\Cal{O}_{K'}} \Cal{O}_{K'_{v_i}}$. As before, the existence of $\Cal{R}$ follows from the local-to-global correspondence. Moreover, by construction, $\Phi_{\Ord}(\Cal{R})=\Cal{R}'$, thus establishing the surjectivity of $\Phi_{\Ord}$. To see that the level ideal of $\Phi_{\Ord}(\Cal{R})$ is equal to $\Phi_{\FI}(\Cal{L}_\Cal{R})$, we simply need to combine a few facts. First, $\Phi_{\Ord}(\Cal{R})$ was defined to have completions everywhere isomorphic to those of $\Cal{R}$. Second, the completion of an order ideal of an $\Cal{O}_K$--module is equal to the order ideal of the completion of the module \cite[Thm 4.20]{Reiner}. In tandem, we obtain the claim on level ideals.
\end{proof}

The following is immediate from Theorem \ref{orderbij} and the definition of the level ideal applied to maximal orders.

\begin{cor}\label{MaxOrderBiject}
If $\mathcal{O}$ is a maximal order in $A$, then $\Phi_{\Ord}(\mathcal O)$ is also a maximal order in $\Phi_{\Br}(A)$.
\end{cor}

\subsection{Arithmetic lattices}

We refer the reader to \cite{Witte} for a general introduction to arithmetic lattices in semisimple Lie groups. 

Given a central simple algebra $A$ over $K$, by the Wedderburn Structure Theorem, $A \cong \Mat(r,D)$, where $D$ is a central simple division algebra. Let $v_{1,\R},\dots,v_{r_1,\R}$ be the real places of $K$ and $v_{1,\C},\dots,v_{r_2,\C}$ be the complex places of $K$, where the latter are taken up to complex conjugation. For each complex place, $A \otimes_K K_{v_{j,\C}} \cong \Mat(rd,\C)$, where $d$ is the degree of $D$ over $K$ while at each real place, we have
\[ A \otimes_K K_{v_{j,\R}} \cong \begin{cases} \Mat(rd,\R),& v_{j,\R} \notin \Ram_\iny(D), \\ \Mat(rd/2,\bf{H}),& v_{j,\R} \in \Ram_\iny(D). \end{cases} \]
The group of norm one elements $A^1$ of $A$ embeds into either $\SL(rd,\C)$, $\SL(rd,\R)$, or $\SL(rd/2,\bf{H})$. Given an order $\Cal{O}$ in $A$, the group of norm one elements $\Cal{O}^1$ embeds into these Lie groups as well. Moreover, by Borel--Harish-Chandra \cite{BHC}, the image of $\Cal{O}^1$ is an arithmetic lattice in the product
\begin{equation}\label{A1}
(\SL(rd,\C))^{r_2} \times (\SL(rd,\R))^{r_1 - \abs{\Ram_\iny(D)}} \times (\SL(rd/2,\mathbf{H}))^{\abs{\Ram_\iny(D)}}.
\end{equation}
Typically, one removes compact factors in the product as the image of $\Cal{O}^1$ is also a lattice in the product of all the non-compact groups. The groups $\SL(rd,\R)$ and $\SL(rd,\C)$ are non-compact provided $rd>1$. The groups $\SL(rd/2,\bf{H})$ are non-compact provided $r>1$ or $d>2$, and so compact only when $rd/2=1$. Additionally, for geometric connections, one typically works with lattices in the adjoint form of (\ref{A1}); we will work with lattices in $A^*/K^*$ below as a result. 
\subsection{Bijections between maximal arithmetic lattices}\label{QuatMaxBijSect}

We restrict our attention to the case in which $A$ is a quaternion algebra and extend the bijection of Theorem \ref{orderbij} to a bijection between the maximal arithmetic subgroups of $A^*/K^*$ and those of $\Phi_{\Br}(A)^*/{K^\prime}^*$. Although this bijection may be obtained by associating to the normalizer $N(\mathcal E)$ of an Eichler order $\mathcal E$ in $A$, the normalizer $N(\Phi_{\Ord}(\mathcal E))$ of the corresponding Eichler order $\Phi_{\Ord}(\mathcal E)$ in $\Phi_{\Br}(A)$ (see \cite[Ch 11.4]{MR} for this characterization of maximal arithmetic subgroups of $A^*/K^*$), it is more natural to work within the context of Bruhat--Tits theory.

We begin by briefly recalling the construction of maximal arithmetic subgroups of quaternion algebras. Our treatment follows that of Borel \cite{borel}, though we will use the somewhat less burdensome notation employed by Chinburg and Friedman \cite[p.~41]{CF}. The Bruhat--Tits tree for $\SL(2,k)$, where $k$ is a non-archimedean local field with ring of integers $\mathcal O_k$ and uniformizer $\pi_k$, is given as follows. Given two maximal orders $\mathcal R_1$ and $\mathcal R_2$ of the split quaternion algebra $\Mat(2,k)$, we define the distance $d(\mathcal R_1,\mathcal R_2)$ to be the non-negative integer $n$ such that as $\mathcal O_k$--modules, $\mathcal R_1/\mathcal R_1\cap\mathcal R_2 \cong \mathcal O_k/\pi_k^n\mathcal O_k$. The vertices of the Bruhat--Tits tree $\mathcal T_k$ for $\SL(2,k)$ are the distinct maximal orders of $\Mat(2,k)$. Two vertices are connected by an undirected edge if the distance between the associated maximal orders is one. We represent edges in the tree $\Cal{T}_k$ by $\{\mathcal E_v,\widehat{\mathcal E}_v\}$. The group $A^*/K^*$ acts on $\mathcal T_{K_v}$ via the conjugation action of $A^*$ on the set of maximal orders of $A_v=A\otimes_K K_{v}$. Let $S$ be a finite set of finite places of $K$ which are disjoint from $\Ram_f(A)$. For a maximal order $\mathcal R$ of $A$, we define 
\[ \Gamma_{\mathcal R, S}:=\set{\overline{x}\in A^*/K^* : x \mbox{ fixes } \mathcal R_v \mbox{ if } v\not\in S, x \mbox{ fixes } \{\mathcal E_v,\widehat{\mathcal E}_v\} \mbox{ when } v\in S}. \]
Borel \cite{borel} has shown that every maximal arithmetic subgroup of $A^*/K^*$ arises in this manner. Denote by $\mathrm{MaxArith}(A,K)$ the set of maximal arithmetic subgroups of $A^*/K^*$.

\begin{prop}\label{MaxBiject}
The bijection $\Phi_{\Ord}$ extends to a bijection between $\mathrm{MaxArith}(A,K)$, $\mathrm{MaxArith}(\Phi_{\Br}(A),K^\prime)$.
\end{prop}

\begin{proof}
Our bijection between $\mathrm{MaxArith}(A,K)$ and $\mathrm{MaxArith}(\Phi_{\Br}(A),K^\prime)$ is the obvious one. Using $\Phi_{\V}$, there is a distance preserving isomorphism between the trees $\mathcal T_{K_v}$ and $\mathcal T_{K^\prime_{\Phi_{\V}(v)}}$ for all $v$ not lying in $\Ram_f(A)$. For $S_\Phi=\{\Phi_{\V}(v): v\in S\}$, the desired bijection sends $\Gamma_{\mathcal R, S}$ to $\Gamma_{\Phi_{\Ord}(\mathcal R), S_{\Phi}}$.
\end{proof}

\section{Volume of the Associated Orbifolds}

In this section, we show that our bijection $\Phi_{\Ord}$ extended to arithmetic lattices derived from maximal orders also preserves covolume. The main tool is a special case of Prasad's volume formula \cite{Prasad} that we work out explicitly.

\subsection{Prasad's volume formula}

We refer the reader to \cite{Prasad} for a thorough treatment of this material. We have also borrowed the notation used in \cite{Prasad} for referencing ease. Let $G$ be an absolutely quasi-simple, simply connected algebraic group defined over a number field $K$ and $\mathcal{G}$ be an absolutely quasi-simple, simply connected algebraic group which is also quasi-split over $K$.  For a place $v\in V^K_f$ let $q_v$ denote the size of the residue field $F_v$, let $S$ be a finite set of places containing $V^K_\infty$, and let $r_v$ be the $k_v^{unr}$--rank of $G$ where $k^{unr}_v$ denotes the maximal unramified extension of $k_v$. We fix a \textbf{coherent} system of parahorics in $G$ by which we mean a collection of parahorics $P_v$, denoted $(P_v)_{v\in V_f^K}$, such that $\prod_{v\in V^K_\infty}G(K_v)\prod_{v\in V^K_f}P_v$ is an open subgroup of the adelic points $G(\mathbf{A}_K)$ (for the definition of parahoric, see \cite{Tits}). Given this coherent system of parahorics, for each $v\in V^K_f$, Bruhat--Tits theory associates a smooth, affine group scheme $G_v$ over $\spec{\mathcal{O}_{K_v}}$ such that the generic fiber $G_v\times_{\mathcal{O}_{K_v}}K_v$ of $G_v$ is isomorphic to the base change of $G$ to $\spec{K_v}$, i.e. $G\times_KK_v$, and further such that the $\mathcal{O}_{K_v}$--points of $G_v$ are isomorphic to $P_v$. Since $G$ is simply connected, the fiber over the closed point in $\spec{\mathcal{O}_{K_v}}$ is connected and we denote it by $\overline{G}_v:=G_v\times_{\mathcal{O}_{K_v}}F_v$. The group $\overline{G}_v$ admits a Levi decomposition over $F_v$ as $\overline{G}_v=\overline{M}_v.R_u(\overline{G}_v)$ and we denote by $\overline{M}_v$ the maximal, connected, reductive part and by $R_u(\overline{G}_v)$ the unipotent radical.  Fixing an $F_v$--defined Borel subgroup $\overline{B}_v$, we let $\overline{T}_v$ denote the maximal $F_v$--torus of $\overline{B}_v$. In the volume formula we will disregard the unipotent radical and only consider the reductive part.  Letters in calligraphy will denote the similar notation for $\mathcal{G}$. The following volume formula of Prasad can be found in \cite[Thm 3.5]{Prasad}.

\begin{thm}[Prasad \cite{Prasad}]\label{PVol}
With the notations above, let $S$ be a finite set of places such that $V^K_\infty\subset S$ and $G_S=\prod_{v\in S}G(K_v)$.
If $\Lambda$ denotes the lattice obtained as the image of $G(K)\cap(G_S\cdot\prod_{v\notin S}P_v)$ under the natural projection to $G_S$, then
\[ \mu_S(G_S/\Lambda)=\mathcal{D}_K^{\frac{1}{2}\dim G}(\mathcal{D}_L/\mathcal{D}_K^{[L:K]})^{\frac{1}{2}s(\mathcal{G})}\left(\prod_{v\in V_\infty}\left|\prod_{i=1}^r\frac{m_i!}{(2\pi)^{m_i+1}}\right|_v\right)\tau_K(G)\mathcal{E} \]
where $|-|_v$ denotes the valuation given by $v$, $\mathcal{E}$ is given by
\[ \mathcal{E}=\prod_{v\in S\cap V^K_f}\frac{q_v^{(r_v+\dim\overline{\mathcal{M}_v})/2}}{|\overline{T}_v(F_v)|}\prod_{v\notin S}\frac{q_v^{(\dim\overline{M_v}+\dim\overline{\mathcal{M}_v})/2}}{|\overline{M}_v(F_v)|}, \]
and with the notations
\begin{enumerate}
\item $\mathcal{D}_K$ (resp. $\mathcal{D}_L$) is the discriminant of $K$ (resp. $L$) and $L$ is the smallest extension of $K$ over which $\mathcal{G}$ splits.
\item $s(\mathcal{G})$ is an integer based on the type of the group, which is $0$ if $\mathcal{G}$ split over $K$.
\item $m_i$ are the exponents of $G$.
\item $\tau_K(G)$ is the Tamagawa number.
\item $\mu_S$ is the product measure on $G_S$.
\end{enumerate}
\end{thm}

\subsection{Volume formula for central simple algebras}

The following is an explicit, special case of Theorem \ref{PVol} for groups arising from central simple algebras.

\begin{prop}\label{P:volumeformula}
Let $D$ be a degree $d$ division algebra, $G=\SL_n(D)$ the norm one elements of the central simple algebra $\Mat(n,D)$,
$\mathcal{R}$ a maximal order in $D$, and $\Lambda=\SL_n(\mathcal{R})$. If $\Ram_f(D)$ denotes the set of finite ramified places of $D$, then the volume of $G/\Lambda$ is given by 
\begin{align*}
\mu(G/\Lambda)=\mathcal{D}_K^{((nd)^2-1)/2}&\left(\prod_{i=1}^{nd-1}\frac{i!}{(2\pi)^{i+1}}\right)^{[k:\Q]}\prod_{i=1}^{nd-1}\zeta_K(i+1)\\
&\cdot \prod_{v\in \Ram_f(D)}\frac{\displaystyle \prod_{i=1}^{nd-1}(q_v^{i+1}-1)}{\displaystyle \prod_{i=1}^{n_v-1}(q_v^{d_v(i+1)}-1)\left(\sum_{i=0}^{d_v-1}q_v^i\right)}.
\end{align*}
Here we use the conventions that $d_v$ is the order of the local invariant $\Inv_v(D)$, $n_vd_v=nd$, and if $n_v=1$ then $\prod_{i=1}^{n_v-1}(q_v^{d_v(i+1)}-1)= 1$.
\end{prop}

\begin{proof}
Our interest is in algebraic groups of type $^1\Aa_{nd-1}$. Namely let $G$ be the $K$--defined algebraic group with group of $E$--points given by $G(E)=\SL_n(D\otimes_K E)$, where $D$ a fixed division algebra defined over $K$ and extension $E/K$.  In this case, the quasi-split $\mathcal{G}$ is given by $\SL_{nd}$. It is well known that $G$ is an inner $K$--form of $\mathcal{G}$ (see \cite[2.2]{PR} for instance). Given a maximal order $\mathcal{R}\subset D$ we consider the $\mathcal{O}_K$--form of $G$ such that $G(\mathcal{O}_K)=\SL_n(\mathcal{R})$ and $G(\mathcal{O}_{K_v})=\SL_n(\mathcal{R}_v)$ for any $v\in V^K$, with the convention that $\mathcal{R}_v=\mathcal{R}\otimes_{\mathcal{O}_{K}}\mathcal{O}_{K_v}$. It is straightforward to check that $P_v=\SL_n(\mathcal{R}_v)$ is a coherent system of parahorics. For an inner form of type $\Aa_{nd-1}$, we know that $L=K$, $\tau_K(G)=1$, and the exponents are given by $m_i=i$ for $1\le i\le nd-1$ (\cite[1.5]{Prasad} and the fact that quasi-split and split are the same in this case). Since we are not considering $S$--arithmetic, i.e. $S=V^K_\infty$, the volume formula greatly simplifies to give 
\[ \mu(G/\Lambda)=\mathcal{D}_K^{1/2((nd)^2-1)}\left(\prod_{i=1}^{nd-1}\frac{i!}{(2\pi)^{i+1}}\right)^{[K:\Q]}\mathcal{E}, \] 
where $\mathcal{E}=\prod_{v\in V^K_f}\frac{q_v^{(\dim\overline{M_v}+\dim\overline{\mathcal{M}_v})/2}}{|\overline{M}_v(F_v)|}$. Let $Q$ denote the set of non-archimedean places for which $\SL_n(D)$ does not split over $K_v$. If $v\notin Q$ then $D$ splits over $K_v$ and $G\cong\mathcal{G}=\SL_{nd}$ over $K_v$. That implies $\overline{M}_v\cong \overline{\mathcal{M}}_v=\SL_{nd}(\mathcal{O}_{K_v})$ and yields the following elementary manipulation
\begin{align*}
\mathcal{E}&=\prod_{v\in V^K_f}\frac{q_v^{(\dim\overline{M_v}+\dim\overline{\mathcal{M}_v})/2}}{|\overline{M}_v(F_v)|} =\prod_{v\in V^K_f\setminus Q}\frac{q_v^{\dim\overline{\mathcal{M}}}}{|\overline{\mathcal{M}}_v(F_v)|}\prod_{v\in Q}\frac{q_v^{(\dim\overline{M_v}+\dim\overline{\mathcal{M}_v})/2}}{|\overline{M}_v(F_v)|}\\
&=\prod_{v\in V^K_f}\frac{q_v^{\dim\overline{\mathcal{M}}}}{|\overline{\mathcal{M}}_v(F_v)|}\prod_{v\in Q}q_v^{(\dim\overline{M_v}-\dim\overline{\mathcal{M}_v})/2}\frac{|\overline{\mathcal{M}}_v(F_v)|}{|\overline{M}_v(F_v)|}
\end{align*}
For convenience sake we write 
\[ \lambda_v=q_v^{(\dim\overline{M_v}-\dim\overline{\mathcal{M}_v})/2}\frac{|\overline{\mathcal{M}}_v(F_v)|}{|\overline{M}_v(F_v)|}, \]
and in the future refer to each $\lambda_v$ as a \textbf{lambda factor}. Since $\mathcal{G}=\SL_{nd}$, we see that $\overline{\mathcal{M}}_v=\SL_{nd}(\mathbf{F}_{q_v})$ and via \cite[Table 1]{Ono} we compute
\begin{align*}
\prod_{v\in V^K_f}\frac{q_v^{\dim\mathcal{M}}}{|\overline{\mathcal{M}}_v(F_v)|}&=\prod_{v\in V^K_f}\frac{\displaystyle q_v^{(nd)^2-1}}{\displaystyle q_v^{(nd-1)(nd)/2}\prod_{i=1}^{nd-1}(q_v^{i+1}-1)}\\
&=\prod_{v\in V^K_f}\frac{\displaystyle 1}{\displaystyle \prod_{i=1}^{nd-1}\pr{1-\frac{1}{q_v^{i+1}}}} =\prod_{i=1}^{nd-1}\zeta_K(i+1)
\end{align*}
If $v\in Q$, then $\mathcal{R}_v$ does not split completely and the index of the local Dynkin diagram is of type $^{d_v}A_{n_vd_v-1}$ according to the classification in \cite[4.3]{Tits}. 

If $n\ge 2$, then the absolute local Dynkin diagram is a cycle of length $nd$ where $\Gal(K_v^{unr}/K_v)$ acts as the cyclic group $\Z/d_v\Z$ by a rotation of the cycle (see \cite[4.3]{Tits}). The relative local Dynkin diagram is hence a cycle of length $n_v$. For the vertex $\mathbf{x}$ corresponding to our parahoric, by Bruhat--Tits theory (\cite[3.5.2]{Tits}), to find the corresponding index of $\overline{M}_v^{ss}$, one deletes the vertices in the orbit of $\mathbf{x}$ under the Galois action as well as all of the edges adjacent to those vertices.  The resulting diagram gives the desired index. Here $\overline{M}_v^{ss}$ denotes the semisimple part of $\overline{M}_v$, which is nothing more than the derived subgroup of $\overline{M}_v$. In particular, in this case $\overline{M}_v^{ss}$ has absolute type 
\[ \underbrace{A_{n_v-1}\times...\times A_{n_v-1}}_{d_v-\mathrm{times}}. \] 
To compute $|\overline{M}_v|$, we write $\overline{M}_v=\overline{M}_v^{ss}.\overline{R}_v(\overline{M}_v)$, where $\overline{R}_v(\overline{M}_v)$ is the radical of $\overline{M}_v$. It can be shown that $\overline{M}_v$ has $F_v$--rank equal to the $K_v$--rank of $G$ which in this case is $nd-1=n_vd_v-1$. Over $K_v$, $\overline{M}_v^{ss}$ has absolute rank $n_vd_v-d_v$ from the product above and hence the radical of $\overline{M}_v$ is a $d_v-1$ dimensional, non-split torus.  Consequently, $R_v(\overline{M}_v)$ is given by the norm torus
\[ \Res^{(1)}_{\mathbf{F}_{q_v^{d_v}}/\mathbf{F}_{q_v}}(\mathbf{G}_m) =\Res_{\mathbf{F}_{q_v^{d_v}}/\mathbf{F}_{q_v}}(\mathbf{G}_m)\cap G. \]
One can compute that 
\[ \abs{\Res^{(1)}_{\mathbf{F}_{q_v^{d_v}}/\mathbf{F}_{q_v}}(\mathbf{G}_m)} = \frac{q_v^{d_v}-1}{q_v-1} = \sum_{i=0}^{d_v-1}q_v^i, \]
which comes from the fact that $\abs{\mathbf{F}_{q_v^{d_v}}^\times} = q_v^{d_v}-1$ and the well known fact that the norm is surjective onto $\mathbf{F}_{q_v}$. One can similarly compute $\abs{\overline{M}_v^{ss}}$ by noting that it has the same order  as the group $\SL_{n_v}(\mathbf{F}_{q_v^{d_v}})$ since 
\[ \SL_{n_v}(\mathcal{R}_v)\otimes_{\mathcal{O}_{K_v}}\mathbf{F}_{q_v}\cong\SL_{n_v}(\mathcal{R}_v\otimes_{\mathcal{O}_{K_v}}\mathbf{F}_{q_v})\cong\SL_{n_v}(\mathbf{F}_{q_v^{d_v}}). \] 
The latter isomorphism is well-known (\cite[1.4]{PR} for instance). Combining \cite[Table 1]{Ono}
 \[ \abs{\overline{M}_v^{ss}}=q_v^{n_vd_v(n_v-1)/2}\prod_{i=1}^{n_v-1}(q_v^{d_v(i+1)}-1) \]
and Lang's isogeny theorem (\cite[p.~290]{PR}), yields 
\begin{align*}
\abs{\overline{M}_v}&= \abs{\overline{M}_v^{ss}}\abs{R(\overline{M}_v)} =\left(q_v^{n_vd_v(n_v-1)/2}\prod_{i=1}^{n_v-1}(q_v^{d_v(i+1)}-1)\right)\pr{\sum_{i=0}^{d_v-1}q_v^i}
\end{align*}
Hence $\dim\overline{M}_v=n_v^2d_v-1$ from which the lambda factors can be computed as
\begin{align*}
\lambda_v&= q_v^{(n_v^2d_v-1-((nd)^2-1))/2}\frac{\displaystyle q_v^{(nd-1)(nd)/2}\prod_{i=1}^{nd-1}(q_v^{i+1}-1)}{\displaystyle q_v^{n_vd_v(n_v-1)/2}\prod_{i=1}^{n_v-1}(q_v^{d_v(i+1)}-1)\left(\sum_{i=0}^{d_v-1}q_v^i\right)}\\
&= q_v^{n_v^2d_v(1-d_v)/2}\frac{\displaystyle q_v^{n_v^2d_v(d_v-1)/2}\prod_{i=1}^{nd-1}(q_v^{i+1}-1)}{\displaystyle \prod_{i=1}^{n-1}(q_v^{d_v(i+1)}-1)\pr{\sum_{i=0}^{d_v-1}q_v^i}}= \frac{\displaystyle \prod_{i=1}^{nd-1}(q_v^{i+1}-1)}{\displaystyle \prod_{i=1}^{n_v-1}(q_v^{d_v(i+1)}-1)\pr{\sum_{i=0}^{d_v-1}q_v^i}}
\end{align*}
Here we are repeatedly using that $n_vd_v=nd$. This completes the proof in the case of $n_v\ge 2$.

Now if $n=1$, a new phenomenon can occur. If $v\notin Q$, we still have $G\cong \mathcal{G}=\SL_d$ over $K_v$ and hence our computation of $\overline{\mathcal{M}}_v$ from above carries through. Similarly if $v\in Q$ and $n_v>1$, the above computation carries through. If $v\in Q$ such that $n_v=1$, then the absolute local Dynkin diagram is still a cycle of length $d$ where the Galois group acts as a cyclic group of order $d$ by a rotation of this cycle.  However, unlike above, the relative local Dynkin diagram is empty since there is only one orbit under this action. Therefore $\SL_1(\mathcal{R}_v)$ is totally anisotropic and so $\overline{M}_v=\Res^{(1)}_{\mathbf{F}_{q_v^{d_v}}/\mathbf{F}_{q_v}}(\mathbf{G}_m)$. This group has order $\sum_{i=0}^{d-1}q_v^i$ and the associated lambda factor is given by:
\begin{align*}
\lambda_v&= q_v^{(d-1-(d^2-1))/2}\frac{\displaystyle q_v^{d(d-1)/2}\prod_{i=1}^{d-1}(q_v^{i+1}-1)}{\displaystyle \sum_{i=0}^{d-1}q_v^i}= q_v^{d(1-d)/2}\frac{\displaystyle q_v^{d(d-1)/2}\prod_{i=1}^{d-1}(q_v^{i+1}-1)}{\displaystyle \sum_{i=0}^{d-1}q_v^i}\\
&= q_v^{d(1-d)/2}\frac{\displaystyle q_v^{d(d-1)/2}\prod_{i=1}^{d-1}(q_v^{i+1}-1)}{\displaystyle \sum_{i=0}^{d-1}q_v^i}= (q_v-1)\prod_{i=1}^{d-2}(q_v^{i+1}-1)= \prod_{i=1}^{d-1}(q_v^{i}-1).
\end{align*}
\end{proof}

We now return to the setting of locally equivalent number fields. As before, let $K,K^\prime$ be a locally equivalent pair, $D$ a division algebra over $K$, and $\mathcal R$ a maximal order of $D$. Let $D^\prime,\mathcal{R}'$ be the associated division algebra and maximal order, respectively, over $K^\prime$ under the correspondences from Theorem \ref{BrauerIsomorphism} and Theorem \ref{orderbij}.

\begin{cor}\label{cor:cor1tovolumeformula}
If $G=\SL_n(D)$, $G'=\SL_n(D')$, $\Lambda=\SL_n(\mathcal{R})$, and $\Lambda'=\SL_n(\mathcal{R}')$, then the volumes of the associated quotients are the same, namely $\mu(G/\Lambda)=\mu(G'/\Lambda')$.
\end{cor}

Indeed these quantities are completely controlled by the local behavior of the number field, division algebra, and maximal order. Additionally, we know that locally equivalent fields share the same discriminant, zeta function, and degree so the result follows. Corollary \ref{cor:cor1tovolumeformula} extends to any order given by Theorem \ref{orderbij} as well.

\begin{cor}
Let $\mathcal{T}$ be any order contained in $\mathcal{R}$ and let $\mathcal{T}'$ be the corresponding order given by the construction of Theorem \ref{orderbij}. If $\Lambda=\SL_n(\mathcal{T})$ and $\Lambda'=\SL_n(\mathcal{T}')$, then $\mu(G/\Lambda)=\mu(G'/\Lambda')$.
\end{cor}

\begin{proof}
This follows immediately from Corollary \ref{cor:cor1tovolumeformula}, as the proof of Theorem \ref{orderbij} makes it clear that the index of $\mathcal T$ in $\mathcal R$ coincides with the index of $\mathcal T^\prime$ in a maximal order $\mathcal R^\prime$ of $D^\prime$.
\end{proof}

We conclude this section with an example showing that our bijection between maximal arithmetic lattices (Proposition \ref{MaxBiject}) does not always preserve covolumes. In essence this is due to the fact that there exist locally equivalent number fields with different class numbers (cf \cite{deSmit-Perlis}). 

\textbf{Example.}
Let $K_1=\Q(\sqrt[8]{799})$ and $K_2=\Q(\sqrt[8]{16\cdot 799})$. It was shown by de Smit and Perlis \cite{deSmit-Perlis} that these two number fields have isomorphic adele rings and different class numbers. Indeed, using Magma \cite{Magma} it is easy to compute that the class number of $K_1$ is $2^{13}$ and the class number of $K_2$ is $2^{14}$. Because $K_1$ and $K_2$ have isomorphic adele rings, their Dedekind zeta functions are equal \cite{Komatsu3}. It is well known that the signature of a number field is determined by the Dedekind zeta function. In this case we see that both $K_1$ and $K_2$ have signature $(2,3)$. For $i=1,2$ let $B_i$ be the unique quaternion division algebra over $K_i$ which is unramified at all finite primes of $K_i$. It is clear that $B_1$ and $B_2$ correspond to one another via the isomorphism in Theorem \ref{BrauerIsomorphism}. Let $\mathcal{O}_1$ be a maximal order of $B_1$ and $\mathcal{O}_2$ be the corresponding (via Theorem  \ref{orderbij}) maximal order of $B_2$. Let $\Gamma_1$ (respectively $\Gamma_2$) denote the image in $\PSL(2,\C)^3$ of $N(\mathcal{O}_1)$ (respectively $N(\mathcal{O}_1)$), where $N(\mathcal{O}_i)$ is the normalizer in $B_i^*$ of $\mathcal{O}_i^*$. Borel \cite{borel} has shown that these are both maximal arithmetic subgroups of $\PSL(2,\C)^3$. The covolumes of these groups are most easily computed using Chinburg and Friedman's \cite[Prop 2.1]{chinburg-smallestorbifold} simplification of Borel's volume formula:
\[ \textrm{CoVolume}(\Gamma_i) = \frac{\mathcal{D}^{\frac{3}{2}}_{K_i}\zeta_{K_i}(2)}{2^{12}\pi^{7}[K_i(B_i):K_i]},\]
where $\mathcal{D}_{K_i}$ is the absolute value of the discriminant of $K_i$, $\zeta_{K_i}(s)$ the Dedekind zeta function of $K_i$, and $K_i(B_i)$ is the maximal extension of $K_i$ which is unramified at all finite primes of $K_i$ and whose Galois group is an elementary abelian group of exponent $2$. It is known that the $[K_i(B_i):K_i]$ coincides with the type number of $B_i$ \cite[p. 37]{CF}. (Recall that the type number of a central simple algebra defined over a number field is the number of isomorphism classes of maximal orders.) Because $\zeta_{K_i}(s)$ determines $\mathcal{D}_{K_i}$, our claim that $\Gamma_1$ and $\Gamma_2$ have different covolumes follows from the fact (easily verified with Magma \cite{Magma}) that the type number of $B_1$ is $128$ and the type number of $B_2$ is $64$.

\section{Greatest Common Divisors and Rigidity}

We now exhibit a few rigidity results regarding natural isomorphisms of Brauer groups. The first rigidity result which we will prove is Theorem \ref{GCDEQUIV} from the introduction. 

\subsection{Proof of Theorem \ref{GCDEQUIV}}

Recall that $K,K'$ are locally GCD equivalent if for every rational prime $p$ which is unramified in $K/\Q$ and $K'/\Q$ we have
\[\gcd(f(v_1/p),\dots,f(v_g/p))=\gcd(f(v'_1/p),\dots,f(v'_{g'}/p)),\]
where $v_1,\dots,v_g$ are the places of $K$ lying above $p$ and $v'_1,\dots,v'_{g'}$ are the places of $K'$ lying above $p$. Here $f(v_i/p)$ (respectively $f(v'_i/p)$) is the inertia degree of $v_i$ (respectively $v'_i$) over $p$.

\begin{proof}[Proof of Theorem \ref{GCDEQUIV}]
We proceed via contradiction assuming $K,K'$ are not locally GCD equivalent. In that case, there is a prime $p_1\in\Z$ which is unramified in $K/\Q$ and $K'/\Q$ such that 
\[ g_{p_1}=\gcd(f(v_1/{p_1}),\dots,f(v_g/{p_1})) \ne \gcd(f(v'_1/{p_1}),\dots,f(v'_{g'}/p_1)) = g'_{p_1}. \] 
Without loss of generality we may assume that $g'_{p_1}<g_{p_1}$. Let $p_2,\dots,p_{g_{p_1}}$ be distinct rational primes which all have the same GCD of local degrees (relative to the extension $K/\Q$) as $p_1$. Set $B$ to be the degree $g_{p_1}$ division algebra defined over $\Q$ whose local invariants are $\frac{1}{g_{p_1}}$ at $p_1,p_2,\dots,p_{g_{p_1}}$ and which is split at all other rational primes. Notice that if $v$ is a place of $K$ which lies above $p_i$ then 
\[\Inv_v(B\otimes_{\Q} K)=[K_v:\Q_{p_i}]\cdot \frac{1}{g_{p_1}}=f(v/p_i)\cdot \frac{1}{g_{p_1}}\in\Z.\] 
By (\ref{ABHN}), we see that $B\otimes_{\Q} K\cong \Mat(g_{p_1},K)$, and hence $[B\otimes_{\Q} K]$ is trivial in $\Br(K)$. Now consider the algebra $\Phi_{\Br}(B\otimes_{\Q}K')=\Res_{K'/\Q}(B)$. As $g'_{p_1}<g_{p_1}$, there is a place $v'$ of $K'$ which lies above $p_1$ for which $f(v'/p_1)$ is not divisible by $g_{p_1}$. Consequently,  
\[ \Inv_{v'}(B\otimes_{\Q}K')=[K'_{v'}:\Q_{p_1}]\cdot \frac{1}{g_{p_1}}=f(v'/p_1)\cdot \frac{1}{g_{p_1}}  \not\in\Z,\] 
and by (\ref{ABHN}), we see that $B\otimes_{\Q}K'$ represents a nontrivial class in $\Br(K')$. As $\Phi_{\Br}(\Res_{K/\Q}([B])) \ne \Res_{K'/\Q}([B])$, we contradict the naturality of $\Phi_{\Br}$.
\end{proof}

\subsection{Proof of Corollary \ref{RIGIDITY}}

Corollary \ref{RIGIDITY} is an immediate consequence of Theorem \ref{GCDEQUIV} and the following proposition.

\begin{prop}\label{gcdgalois}
Let $K$ and $K'$ be number fields which are locally GCD equivalent. If $K'/\Q$ is Galois then $K'\subset \widehat{K}$ where $\widehat{K}$ is the Galois closure of $K$ over $\Q$.
\end{prop}

\begin{proof}
Given a rational prime $p$ which is unramified in $K'/\Q$ and splits completely in $\widehat{K}/\Q$, we have $f(v/p)=1$ for all places $v$ lying over $p$. As $K,K'$ are locally GCD equivalent, we see that
\[ 1=f(v'_1/p)=\cdots=f(v'_{g'}/p), \] 
where $v_1,\dots,v_{g'}$ are the distinct places of $K'$ lying over $p$. Thus $p$ splits completely in $K'/\Q$, and so all but finitely many primes of $\Q$ which split completely in $\widehat{K}/\Q$ also split completely in $K'/\Q$. The proof is finished with a standard consequence of the Chebotarev density theorem (cf. \cite[Thm 9, p.~168]{Lang-ANT}).
\end{proof}

Our final result of this subsection is the following rigidity result.

\begin{thm}\label{GaloisBrauerRigidity2}
Let $K, K'$ be number fields and $B, B'$ be central division algebras over $K, K'$ of degree $d$ such that for every field $F\subset K\cap K'$ and division algebra $B_0$ over $F$, $B_0\otimes_F K = B$ if and only if $B_0\otimes_F K' = B'$. If there exists a common subfield $F\subset K\cap K'$ such that $K/F$ and $K'/F$ are both Galois of degree dividing $d$ and $\Res_{K/F}^{-1}(B)\neq \emptyset$ then $K=K'$ and $B=B'$.
\end{thm}

\begin{proof}
By hypothesis there exists a subfield $F\subset K\cap K'$ such that $\Res_{K/F}^{-1}(B)\neq \emptyset$. We assume $F=\bf{Q}$ for simplicity as the general case can be argued identically. To prove that $K=K'$, we will show that rational primes have the same splitting behavior over $K,K'$ using central simple algebras. From $K=K'$, it is a simple matter to deduce $B=B'$. We now commence with the proof.

Given $\widetilde{B}\in \Res_{K/\bf{Q}}^{-1}(B)$, we select a rational prime $p_0$ that is unramified in both $K/\bf{Q}$ and $K'/\bf{Q}$, does not lie below a place of $K$ or $K'$ which ramifies in $B$ or $B'$, does not ramify in $\widetilde{B}$, and does not split completely in $K/\bf{Q}$. As $p_0$ neither ramifies nor splits completely in $K/\bf{Q}$, every place $v$ of $K$ lying above $p_0$ has inertial degree $f$ for some $f>1$. If $g$ is the number of places of $K$ lying above $p_0$, then $fg=[K:\bf{Q}]$. We now select $f-1$ additional primes $p_1,\dots, p_{f-1}$ under identical constraints. We note that the existence of these primes follows from the Chebotarev density theorem. By (\ref{ABHN}), there exists a division algebra $B_0$ over $\bf{Q}$ whose local invariants coincide with those of $\widetilde{B}$ at the primes of $\bf{Z}$ which ramify in $\widetilde{B}$ and which has local invariant $\frac{1}{f}$ at the primes $p_0,\dots, p_{f-1}$. If $v$ is a place of $K$ which lies above one of the $p_i$, then  $\Inv_v(B_0\otimes_{\bf{Q}} K)=1$. It follows that $B_0\otimes_{\textbf{Q}} K=B$, and so by hypothesis, we must also have $B_0\otimes_{\textbf{Q}} K'=B'$. We assert that $p_0$ does not split completely over $K'$. Assuming the contrary, for any place $v'$ of $K'$ which lies above $p_0$, we see that 
\[ B'\otimes_{K'} K'_{v'}=B_0\otimes_{\bf{Q}} K'\otimes_{K'} K'_{v'}=B_0\otimes_{\bf{Q}} K'_{v'}=B_0\otimes_{\bf{Q}}\textbf{Q}_{p_0}. \] 
Since $B_0$ was defined to have local invariant $\frac{1}{f}$ at $p_0$, we conclude that $B_0\otimes_{\bf{Q}}\textbf{Q}_{p_0}$ is a division algebra. In particular, $B'\otimes_{K'} K'_{v'}$ represents a nontrivial class in $\Br(K'_{v'})$ and so $v'$ ramifies in $B'$. However, by selection, $p_0$ does not lie below any prime in $K'$ that resides in $\Ram(B')$. Having obtained a contradiction, we see that $p_0$ does not split completely over $K'$. In total, with a finite number of exceptions, if a rational prime does not split completely in $K/\bf{Q}$ then it does not split completely in $K'/\bf{Q}$. Equivalently, if a rational prime splits completely in $K'/\bf{Q}$ then it splits completely in $K/\bf{Q}$ (with at most a finite number of exceptions). Via the same argument, with the roles of $K$ and $K'$ interchanged, we see that the set of rational primes splitting completely in $K/\bf{Q}$ coincides with the set of rational primes splitting completely in $K'/\bf{Q}$ with at most finitely many exceptions. The Chebotarev density theorem then implies that $K=K'$. Finally, as $B = \widetilde{B} \otimes_\Q K = \widetilde{B} \otimes_\Q K' = B'$, we see that $B=B'$.
\end{proof}

A similar result was proven in the context of quaternion algebras defined over number fields with a unique complex place in \cite[Thm 1.1]{McReid}. Theorem \ref{GaloisBrauerRigidity2} generalizes that result to division algebras of arbitrary degree.

\section{Proof of Theorem \ref{EffectivePR}}

We first note that the work of Prasad and Rapinchuk \cite{PrasadRap} shows that it suffices to obtain an upper bound on the number of isomorphism classes of division algebras defined over $K$ which possess precisely the same set of maximal subfields as $D$. In particular, non-commensurable, length commensurable manifolds must arise from division algebras defined over $K$ (i.e. the associated arithmetic lattices are commensurable with $\SL_1(\mathcal{R})$ for a maximal order $\mathcal{R}$ in a division algebra $D'$ over $K$).

For a division algebra $D$ over $K$, we set $\mathrm{Gen}(D)$ to be the number of isomorphism classes of division algebras over $K$ with the same maximal subfields as $D$. Let $\Ram_f(D)=\{v_1,...,v_n\}$ be the set of finite places of $K$ which ramify in $D$.  We write $d_{v_i}=\deg(D_{v_i})$ and call $d_{v_i}$ the \textbf{local degree} at $v_i$. We set $\Theta_D = \prod_{i=1}^n \phi(d_{v_i})$, where $\phi$ is the Euler $\phi$--function, and note that $\mathrm{Gen}(D) \leq \Theta_D$. Observe that $\prod_{i=1}^{d-1} \zeta_K(i+1)\geq 1$, $\prod_{i=1}^{d-1}\frac{i!}{(2\pi)^{i+1}}>1$ for $d>28$. Moreover, it is straightforward to see that the second product is always greater than $10^{-33}$. It follows from Proposition \ref{P:volumeformula} that
\begin{equation}\label{E:LambdaVol}
\prod_{v\in\Ram_f(D)}\lambda_v\leq \frac{10^{33} V}{\mathcal{D}_K^{(d^2-1)/2}}\leq 10^{33}V.
\end{equation}
Theorem \ref{EffectivePR} is a direct consequence of (\ref{E:LambdaVol}) and the following proposition.

\begin{prop}\label{isobound}
For a fixed $N\in \N$, let $\alpha \in \N$ be such that $3^\alpha \leq N < 3^{\alpha + 1}$ and let $D$ be a division algebra over $K$ such that $\prod_{v\in \Ram_f(D)}\lambda_v\le N$. 
\begin{itemize}
\item[(a)]
If $\alpha < 2$, then $\mathrm{Gen}(D)=\Theta_D = 1$.
\item[(b)]
If $\alpha \geq 2$, then $\Theta_D \le 2^\alpha$.
\end{itemize}
\end{prop}

\begin{proof} 
For (a), observe that the only division algebras $D$ satisfying $\prod_{v\in \Ram_f(D)}\lambda_v\le N\le 8$ must have local degree $2$ at all finite ramified places. Thus, $D$ is a quaternion algebra and it follows that $\mathrm{Gen}(D)=\Theta_D = 1$.

For (b), we will find a maximizer of $\Theta_D$ subject to $\prod_{v\in \Ram_f(D)}\lambda_v\le N$ as we vary over $K$ and $D$. First we will reduce to analyzing a number field $K$ with $[K:\Q]\ge\alpha$ in which the rational prime $2$ splits completely and to a division algebra $D$ ramified only at the places lying above $2$. More specifically we want to reduce to the case where $q_v=2$ for each $v\in \Ram_f(D)$.  Recall that
\[ \lambda_v= \begin{cases} \displaystyle \prod_{i=1}^{d-1}(q_v^i-1), n_v=1\\ \displaystyle \prod_{i=1, d_v\nmid i}^{d-1}(q_v^i-1), n_v\neq 1. \end{cases} \]
For fixed $d_v$, $\lambda_v$ is visibly smallest when $q_v=2$. Furthermore, by definition of $\Theta_D$, changing the size of $q_v$ does not change the value of $\Theta_D$. Hence to minimize the possible values of $\lambda_v$ while simultaneously maximizing $\Theta_D$, we assume $q_v=2$ for all $v$ and only ramify the division algebra at places above $2$ (that does not necessarily force 2 to split completely but we may as well assume it does). The requirement $[K:\Q]\ge \alpha$ is evident from the above discussion and the specific relationship between $\alpha$ and $N$ (see the statement of Proposition \ref{isobound}).

Subject to the above reductions, we now deduce the upper bound by explicitly constructing a division algebra which maximizes $\Theta_D$. That will be accomplished by finding a cubic, division algebra $D_{max}$ that ramifies at as many places lying above the prime $2$ as possible. We assume now that $K$ is a number field such that $[K:\Q]\ge\alpha$, such that $2$ splits completely in $K$, and that $D_{max}$ is a cubic, division algebra with $\alpha$ of the places above $2$ having local degrees $d_{v_i}=3$. By choosing some of the local invariants $\Inv_v(D_{max})$ to be $2/3$ instead of $1/3$, we can always ensure the existence of $D_{max}$ by (\ref{ABHN}),  so long as $\alpha \geq 2$. By construction of $D_{max}$ and by definition of $\alpha$, we see that $\prod_{v\in \Ram_f(D_{max})}\lambda_v\le N$. We claim that $\Theta_{D_{max}}$ is the maximum value of $\Theta_D$ subject to the constraint on lambda factors.

Consider another algebra $D$ such that $\prod_{v\in\Ram_f(D)}\lambda_v\le N$, which we may also assume to be only ramified at primes above $2$. If the degree of $D$ is $3$ then $\Theta_D\le \Theta_{D_{max}}$ by construction. For algebras $D$ of higher degree, we will reduce to the cubic case by finding a cubic, division algebra $D'$ such that $\prod_{v \in \Ram_f(D')} \lambda_v \le N$ and consequently $\Theta_D \le \Theta_{D'} \le \theta_{D_{max}}$. To that end, we let $\lambda_{v}$ be a lambda factor for a place $v\in \Ram_f(D)$ and construct a cubic, division algebra $D'$ as follows:

\begin{enumerate}
\item[(1)] For any place $v\in \Ram_f(D)$ such that $d_v=2$ we will not ramify $D'$, as these do not affect $\Theta_D$.
\item[(2)] For any $v\in \Ram_f(D)$ such that $d_v=3$, we ramify $D'$ at the same place with $d_v=3$.  As remarked above, one can always ensure (\ref{ABHN}) while having all local degrees equal to $3$.
\item[(3)] For any place $v\in \Ram_f(D)$ such that $d_{v}>3$ (necessarily at least $2$ such exist) we ramify $D'$ at as many places over the prime $2$ as possible with local degree $3$ and such that the product of the lambda factors of these primes remains less than $\lambda_{v}$.
\end{enumerate}

By construction, $D'$ is a cubic, division algebra with $\prod_{v\in \Ram_f(D)}\lambda_v\ge \prod_{v\in \Ram_f(D')}\lambda_v$ and so only the claim on $\Theta_D$--values needs to be verified.

\textbf{Claim 1.}
$\Theta_D\le\Theta_{D'}$. 

\begin{proof}[Proof of Claim 1]
To prove the claim, we will call a place $v\in\Ram_f(D)$ \textbf{type (1), (2), or (3)} in correspondence with which item it falls into in the above list. If $v\in \Ram_f(D)$ is of type (3), then
\begin{align*}
\lambda_{v}
&=\prod_{i=1, d_{v}\nmid i}^{d-1}(2^i-1)=\prod_{i=1}^{2}(2^i-1)\prod_{i=3, d_{v}\nmid i}^{d-1}(2^i-1)= 3\prod_{i=3, d_{v}\nmid i}^{d-1}(2^i-1).
\end{align*}
Note that $\prod_{i=1}^{2}(2^i-1)$ is precisely the lambda factor for a cubic, division algebra and a place lying over $2$ for which $d_v=3$. For $i>4$ we know that $2^i-1\ge 2^{i-1}\ge3^{\textstyle \ceil{\frac{i-1}{2}}}$. Using the convention that, if $5>d-1$, $\prod_{i=5, d_{v}\nmid i}^{d-1}\ceil{\frac{i-1}{2}} = 0$, we see that
\[ \lambda_{v}=3\prod_{i=3, d_{v}\nmid i}^{d-1}(2^{i}-1)\ge3^{\textstyle \pr{2+c_4+\sum_{i=5, d_{v}\nmid i}^{d-1}\ceil{\frac{i-1}{2}}}},\]
where $c_4$ is 0 if $d_v\mid 4$ and 2 otherwise. In the above product, $c_4$ arises from finding the biggest power of $3$ less than $2^i-1$ for $i=4$. Using this inequality we see that for each $v\in\Ram_f(D)$ with $d_v>3$, our procedure will instead ramify at least $2+c_4+\sum_{i=5, d_{v}\nmid i}^{d-1}\ceil{\frac{i-1}{2}}$  
places $v'$ above $2$, all with local degrees $d_{v'}=3$.

We now verify $\Theta_D\le \Theta_{D'}$. If a given $v\in\Ram_f(D)$ is of type (1) or (2), then the contribution to $\Theta_D$ from $v$ is the same as the corresponding $v'\in\Ram_f(D')$ to $\Theta_{D'}$. Thus it suffices to deal with $v\in\Ram_f(D)$ of type (3). First note trivially that $\phi(d_v)\le d_v-1\le d-1$. If $v\in\Ram_f(D)$ is of type (3), then the contribution to $\Theta_{D'}$ from the multiple lambda factors of the corresponding $v'\in\Ram_f(D')$ is at least 
\[ \phi(3)^{\textstyle\pr{2+c_4+\sum_{i=5, d_{v}\nmid i}^{d-1}\ceil{\frac{i-1}{2}\}}}}=2^{\textstyle \pr{2+c_4+\sum_{i=5, d_{v}\nmid i}^{d-1}\ceil{\frac{i-1}{2}}}}. \]
So long as this quantity is greater than $\phi(d_v)$ for each $v\in\Ram_f(D)$ of type (3), then we complete the proof.
Indeed if $d_v=4$, then we are done since by our convention 
\[ 2^{\textstyle \pr{2+c_4+\sum_{i=5, d_{v}\nmid i}^{d-1}\ceil{\frac{i-1}{2}}}}=2^{\textstyle 2+c_4}\ge 3=d_v-1. \]
We must have $d_v\mid d$ and consequently $d_v\nmid d-1$ for any division algebra. Hence if $d_v\ge 5$, we again conclude
\[ 2^{\textstyle \pr{2+c_4+\sum_{i=5, d_{v}\nmid i}^{d-1}\ceil{\frac{i-1}{2}}}}\ge 2^{\textstyle \pr{\ceil{\frac{d+2}{2}}}}\ge d_v-1. \]
\end{proof}
\end{proof}

\section{Galois Cohomological Bijections}

We extend the bijection in Theorem \ref{BrauerIsomorphism} to various Galois cohomology sets as well as maximal arithmetic lattices in inner forms of absolutely almost simple, $\Q$--split algebraic groups. 

\subsection{Galois cohomology and forms of algebraic groups}

We denote the $i$th Galois cohomology set with coefficients in the algebraic group $G$ by $H^i(K,G)=H^i(\Gal(\overline{K}/K),G)$ with the understanding that when $G$ is not abelian we will only take $i\in\{0,1\}$. For a number field $K$ and an absolutely almost simple, $K$--split algebraic group $G$, we have a map $H^1(K,\overline{G}) \to H^1(K,\Aut_{\overline{K}}(G))$ where $\overline{G}$ denotes the corresponding adjoint group. Twisting $G$ by a class in $H^1(K,\Aut_{\overline{K}}(G))$ gives a $K$--form of $G$ that is inner when the class is in the image of this map. In this section, we construct a natural bijection between $H^1(K,\overline{G})$ and $H^1(K',\overline{G})$ which, after identification, gives the requisite bijection between inner $K$--forms of $G$ with inner $K'$--forms of $G$. 

It is well known that there is a group isomorphism between $H^2(K,\overline{K}^*)$ and $\Br(K)$, and consequently a group isomorphism $H^2(K,\mu_n)\cong \Br_n(K)$, where $\Br_n(K)$ denotes the $n$--torsion in the Brauer group (see \cite[p.~351]{NSW} or \cite[p.~28]{PR} for instance). Hence, if $K,K'$ are a locally equivalent pair, then there is an isomorphism $\Phi_{\mathrm{GC}}\colon H^2(K,\overline{K}^*)\to H^2(K',\overline{K'}^*)$ induced by the bijection of places $\Phi_V\colon V^K\to V^{K'}$ and such that
\[ \xymatrixrowsep{.2in} \xymatrixcolsep{.2in}\xymatrix{H^2(K,\overline{K}^*) \ar[r]^-{\Phi_{\mathrm{GC}}}\ar[d]_{\iota_K} & H^2(K',\overline{K'}^*)\ar[d]^{\iota_{K'}}\\ \prod_vH^2(K_v,\overline{K}^*) \ar[r]^-{\prod_v\Phi_v}& \prod_{v'}H^2(K'_{v'},\overline{K'}^*)} \]
commutes. Here the vertical arrows are the canonical maps into the places and the $\prod_v\Phi_v$ are isomorphisms induced by the $\Phi_V$.
Furthermore, $\Phi_{\mathrm{GC}}$ is natural with respect to changing fields. We also have an isomorphism $\Phi_{\mathrm{GC},n}\colon H^2(K,\mu_n)\to H^2(K',\mu_n)$ induced by $\Phi_V\colon V^K\to V^{K'}$ such that the diagram
\[ \xymatrix{H^2(K,\mu_n) \ar[r]^-{\Phi_{\mathrm{GC},n}}\ar[d]_{\iota_{\mu_n,K}}& H^2(K',\mu_n)\ar[d]^-{\iota_{\mu_n,K'}}\\
 \prod_vH^2(K_v,\mu_n) \ar[r]^-{\prod_v\Phi_v}& \prod_{v'}H^2(K',\mu_n)} \]
commutes. Again the vertical arrows are the canonical maps into the places and the $\prod_v\Phi_v$ are isomorphisms induced by the $\Phi_V$.

\subsection{Proof of Theorem \ref{GC-Isomorphisms}}

Throughout this section we assume that $G$ is an absolutely almost simple, $\Q$--split algebraic group. Given this, we have the exact sequence
\begin{equation}\label{E:ExactSeq}
\xymatrix{1\ar[r]& Z\ar[r]&\widetilde{G}\ar[r]&\overline{G} \ar[r]&1},
\end{equation}
where $\overline{G}$ denotes the adjoint form of $G$, $\widetilde{G}$ denotes the simply connected form, and $Z$ denotes the fundamental group of $\overline{G}$. We have the associated cohomology exact sequence
\[ \xymatrixrowsep{.2in} \xymatrixcolsep{.2in} \xymatrix{H^1(K,\widetilde{G})\ar[r]\ar[d] &H^1(K,\overline{G})\ar[r]^\delta \ar[d]&H^2(K,Z)\ar[d]  \\ \prod_vH^1(K_v,\widetilde{G})\ar[r]  &\prod_vH^1(K_v,\overline{G})\ar[r]&\prod_vH^2(K_v,Z).} \]

\begin{lemma}\label{LocalEquivCoho}
If $K,K'$ are a locally equivalent pair and $v'=\Phi_V(v)$, then the local diagram
\[ \xymatrixrowsep{.2in} \xymatrixcolsep{.2in}\xymatrix{H^1(K_v,\widetilde{G})\ar[d]_-{\varphi_{v,0}}\ar[r]^-{\gamma_v}&H^1(K_v,\overline{G})\ar[d]_{\varphi_{v,1}}\ar[r]^-{\delta_v}&H^2(K_v,Z)\ar[d]^-{\varphi_{v,2}} \\
H^1(K'_{v'},\widetilde{G})\ar[r]^-{\gamma'_{v'}}&H^1(K'_{v'},\overline{G})\ar[r]^{\delta'_{v'}}&H^2(K'_{v'},Z)~,} \]
commutes, where the vertical maps are bijections induced by local equivalence.
\end{lemma}

\begin{proof}
The isomorphism $K_v\cong K'_{v'}$ induces an isomorphism between $\Gal(\overline{K}_v/K_v)$ and $\Gal(\overline{K'}_{v'}/K'_{v'})$ giving the vertical arrows. As Galois cohomology is functorial with respect to change of group maps, we obtain the commutativity of the diagram (\cite[III]{Berhuy} for instance). 
\end{proof}

Combining the above diagrams, we obtain the following commutative Galois cohomological diagram:
\begin{equation}\label{E:BigGC-Diagram}
\xymatrixrowsep{.2in} \xymatrixcolsep{.2in} \xymatrix{
H^1(K,\widetilde{G})\ar[d]_-{\iota_{\widetilde{G},K}}\ar[r]^-{\gamma}&H^1(K,\overline{G})\ar[r]^-{\delta}\ar[d]_{\iota_{\overline{G},K}}&H^2(K,Z)\ar[d]^-{\iota_{Z,K}}\\ \prod_vH^1(K_v,\widetilde{G})\ar[r]^-{\gamma_V}\ar[d]_-{\varphi_0}& \prod_{v}H^1(K_v,\overline{G})\ar[d]_{\varphi_1}\ar[r]^-{\delta_V}&\prod_{v}H^2(K_v,Z)\ar[d]^-{\varphi_2} \\ \prod_{v'}H^1(K'_{v'},\widetilde{G})\ar[r]^-{\gamma'_{V'}}& \prod_{v'}H^1(K'_{v'},\overline{G})\ar[r]^-{\delta'_{V'}}&\prod_{v'}H^2(K'_{v'},Z) \\ H^1(K',\widetilde{G})\ar[u]^-{\iota_{\widetilde{G},K'}}\ar[r]^-{\gamma'}&H^1(K',\overline{G})\ar[r]^-{\delta'}\ar[u]^{\iota_{\overline{G},K'}}&H^2(K',Z)\ar[u]_-{\iota_{Z,K'}}.}
\end{equation}
In the notation of Lemma \ref{LocalEquivCoho}, we write $\gamma_V=\prod_v\gamma_v$, $\delta_V=\prod_v\delta_v$, and $\varphi_i=\prod_v\varphi_{v,i}$ for $i\in\{0,1,2\}$ with similar notation for the other maps. The following is the main result of this subsection.

\begin{thm}\label{AdBijectionThm}
If $K,K'$ are a locally equivalent pair, then there is a bijection $\Phi_{\Ad}:H^1(K,\overline{G})\to H^1(K',\overline{G})$ such that for $\xi\in H^1(K,\overline{G})$ and $\xi'=\Phi_{\Ad}(\xi)\in H^1(K',\overline{G})$, the image $(\xi_v)$ of $\xi$ under $\iota_{\overline{G},K}$ is mapped by $\varphi_1$ to the image $(\xi'_{v'})$ of $\xi'$ under $\iota_{\overline{G},K'}$.
\end{thm}

We briefly remark that, as in the proofs in Section \ref{BrauerIsoSection}, the proof of Theorem \ref{AdBijectionThm} depends on the initial choice of $\Phi_V$ and the isomorphisms $K_v\cong K'_{\Phi_V(v)}$. 

\begin{proof}
Before beginning the proof in full, we briefly outline our strategy. Since $\overline{G}$ satisfies the Hasse principle, the maps $\iota_{\overline{G},K},\iota_{\overline{G},K'}$ are injective. In particular, if
\begin{equation}\label{E:GC-Image}
\vp_1\pr{\iota_{\overline{G},K}\pr{H^1(K,\overline{G})}} \subset \iota_{\overline{G},K'}\pr{H^1(K',\overline{G})},
\end{equation}
then we can set $\Phi_{\Ad} = \iota_{\overline{G},K'}^{-1} \circ \vp_1 \circ \iota_{\overline{G},K}$. To establish (\ref{E:GC-Image}), for a class $\xi \in H^1(K,\overline{G})$, we must find the unique class $\xi' \in H^1(K',\overline{G})$ with $\varphi_1(\iota_{\overline{G},K}(\xi)) = \iota_{\overline{G},K'}(\xi')$. To find such a class, we use natural bijections $\Phi_{\widetilde{\Ad}}$ and $\Psi$ between $H^1(K,\widetilde{G})$, $H^1(K',\widetilde{G})$ and $H^2(K,Z)$, $H^2(K',Z)$, respectively. In the event that $\delta(\xi)$ is the trivial class in $H^2(K,Z)$, by exactness we can lift $\xi$ to a class $\widetilde{\xi} = \gamma^{-1}(\xi) \in H^1(K,\widetilde{G})$. We then apply $\Phi_{\widetilde{\Ad}}$ and $\gamma'$ to obtain $\xi' = \gamma'(\Phi_{\widetilde{\Ad}}(\widetilde{\xi}))$. In the case $\delta(\xi)$ is not the trivial class, we twist by a certain cocycle to arrange for the image $\xi$ under the twisted counterpart of $\delta$ to have trivial image in $H^2(K,Z)$. Paired with the natural maps between the twisted/untwisted cohomology sets, we find $\xi'$ in the general case. Interchanging the roles of $K,K'$ in the above argument yields the reverse containment for (\ref{E:GC-Image}). With our outline complete, we now commence with the proof. To begin, we have the following portion of the Galois cohomological diagram (\ref{E:BigGC-Diagram}):
\begin{equation}\label{E:MediumGC-Diagram} 
\xymatrixrowsep{.2in} \xymatrixcolsep{.2in} \xymatrix{ H^1(K,\overline{G})\ar[d]_{\iota_{\overline{G},K}}\ar[r]^-{\delta}&H^2(K,Z)\ar[r]\ar[d]^{\iota_{Z,K}}&1\\ \prod_vH^1(K_v,\overline{G})\ar[r]^-{\delta_V}\ar[d]_{\varphi_1}&\prod_vH^2(K_v,Z)\ar[d]^-{\varphi_2}  \\ \prod_{v'}H^1(K'_{v'},\overline{G})\ar[r]^-{\delta'_{V'}}&\prod_{v'}H^2(K'_{v'},Z) \\ H^1(K',\overline{G})\ar[u]^{\iota_{\overline{G},K'}}\ar[r]^-{\delta'}&H^2(K',Z)\ar[r]\ar[u]_{\iota_{Z,K'}}&1}
\end{equation}
As $\overline{G}$ satisfies the Hasse principle (see \cite[Thm 6.22]{PR}), the maps $\iota_{\overline{G},K}$ and $\iota_{\overline{G},K'}$ are injective.
The maps $\delta, \delta'$ are surjective (see \cite[Thm 6.20]{PR}), while Lemma \ref{LocalEquivCoho} implies that $\varphi_1, \varphi_2$ are bijective. We now define the following classes:
\begin{align}\label{E:Define-Classes}
\eta=\delta(&\xi), \quad  (\xi_v) = \iota_{\overline{G},K}(\xi), \quad (\xi'_{v'})=\varphi_1((\xi_v)) = (\varphi_{v,1}(\xi_v)), \\
&(\eta_v)=\iota_{Z,K}(\eta), \quad (\eta'_{v'})=\varphi_2((\eta_v))=(\varphi_{v,2}(\eta_v)).  \notag
\end{align}

We first construct a natural bijection between $H^2(K,Z)$ and $H^2(K',Z)$. As $G$ is absolutely almost simple and inner, we know that $Z=\mu_n$ for some $n$ or $Z=\mu_2\times\mu_2$ in the case of type $\D_{2k}$ (see \cite[p.~332]{PR}). For any type other than $\D_{2k}$, there is a natural bijection $\Phi_{\mathrm{GC},n}\colon H^2(K,Z)\to H^2(K',Z)$ for such $Z$. Thus there exists a corresponding $\eta'=\Phi_{\mathrm{GC},n}(\eta)$ such that $\eta'$ has image $(\eta'_{v'})$ under $\iota_{Z,K'}$. For type $\D_{2k}$, we have $Z=\mu_2\times\mu_2$, and so there is a natural isomorphism $H^2(K,Z)\cong H^2(K,\mu_2)\times H^2(K,\mu_2)$. Furthermore that isomorphism is functorial in the sense that these maps also commute with the change of group maps. Additionally, the maps $\iota_{Z,K}$ and $\iota_{Z,K'}$ are injective for the above $Z$ by class field theory. Hence we induce a bijection $\Psi\colon H^2(K,Z)\to H^2(K',Z)$ enjoying the same naturality property. Specifically, if $\eta\in H^2(K,Z)$ and $\eta'=\Psi(\eta)$, then $(\varphi_2\circ \iota_{Z,K})(\eta)=\iota_{Z,K'}(\eta')$. The following is the main step in the construction of $\Phi_{\Ad}$.

\textbf{Claim 2.} 
\emph{There exists a unique $\xi'\in(\delta')^{-1}(\eta')$ such that $\iota_{\overline{G},K'}(\xi') = \varphi_1(\iota_{\overline{G},K}(\xi)) = (\xi'_{v'})$.}

Momentarily assuming Claim 2, we prove of Theorem \ref{AdBijectionThm}. By Claim 2, there exists a unique $\xi'\in(\delta')^{-1}(\eta')$ such that $\iota_{\overline{G},K'}(\xi')=(\xi'_v)$. We define $\Phi_{\mathrm{GC}}(\xi) = \xi'$.  As $\iota_{\overline{G},K'}$ is injective, $\Phi_{\mathrm{GC}}$ is injective. Interchanging the roles of $K,K'$ we see that $\Phi_{\mathrm{GC}}$ is surjective, completing our proof of Theorem \ref{AdBijectionThm}. 

\begin{proof}[Proof of Claim 2]
From the commutativity of the bottom square of diagram (\ref{E:MediumGC-Diagram}), there is $\zeta'\in H^1(K',\overline{G})$ such that $\delta'_{V'}(\iota_{\overline{G},K'}(\zeta'))=\delta'_{V'}((\xi'_{v'}))$. Twisting the exact sequence (\ref{E:ExactSeq}) by $\zeta'$, we obtain the sequence
\[ \xymatrix{1\ar[r]&Z\ar[r]&_{\zeta'}\widetilde{G}\ar[r]&_{\zeta'}\overline{G}\ar[r]&1}. \]
The twisted version of (\ref{E:BigGC-Diagram})  is given below with the associated twisted maps decorated with $\zeta'$:
\begin{equation}\label{E:MediumGC-Diagram2}
\xymatrixrowsep{.2in} \xymatrixcolsep{.2in} \xymatrix{
H^1(K',_{\zeta'}\widetilde{G})\ar[d]_-{\iota_{\widetilde{G},K',\zeta'}}\ar[r]^-{\gamma'_{\zeta'}}&H^1(K',_{\zeta'}\overline{G})\ar[r]^-{\delta'_{\zeta'}}\ar[d]_{\iota_{\overline{G},K',\zeta'}}&H^2(K',Z)\ar[d]^-{\iota_{Z,K',\zeta'}} \ar[r] & 1 \\
\prod_{v'}H^1(K'_{v'},_{\zeta'}\widetilde{G})\ar[r]^-{\gamma'_{V',\zeta'}}& \prod_{v'}H^1(K'_{v'},_{\zeta'}\overline{G})\ar[r]^-{\delta'_{V',\zeta'}}&\prod_{v'}H^2(K'_{v'},Z).}
\end{equation}
We have a natural bijection $\tau_{\zeta'}\colon H^1(K',_{\zeta'}\overline{G})\to H^1(K',\overline{G})$ which takes the class of the trivial cocycle to $\zeta'$ (\cite[Prop 35]{Serre} for instance). We also have a map $\tau_{\eta'}\colon H^2(K',Z) \to H^2(K',Z)$ given by multiplication by the class of $\eta'$. These maps are functorial with respect to the connecting map and its twisted counterpart (\cite[Prop II.5.6]{Berhuy}). In total we obtain the following diagram:
\begin{equation}\label{E:TheCube}
\xymatrixrowsep{.45in}
\xymatrixcolsep{.6in}
\xymatrix@!0{
& H^1(K,\widetilde{G}) \ar[dd]_{\Phi_{\widetilde{\Ad}}} \ar[rr] & & H^1(K,\overline{G}) \ar[rr]^{\delta}\ar@{-->}[dd]_{\Phi_{\Ad}} & & H^2(K,Z)\ar[dd]^{\Psi} \ar[rr] & & 1 \\ & & & & & \\
& H^1(K',\widetilde{G}) \ar[rr] & & H^1(K',\overline{G}) \ar[rr]^{\delta'}\ar'[d][dd] & & H^2(K',Z)\ar'[d][dd]^{\iota_{Z,K'}}\ar[rr] & & 1 \\ H^1(K',_{\zeta'}\widetilde{G})\ar[rr]^{\gamma'_{\zeta'}} \ar[dd]_{\iota_{\widetilde{G},K',\zeta'}} & &H^1(K',_{\zeta'}\overline{G}) \ar[ur]^-{\tau_{\zeta'}}\ar[rr]\ar[dd]_{\iota_{\overline{G},K',\zeta'}} & & H^2(K',Z)\ar[ur]_-{\tau_{\eta'}}\ar[dd] \ar[rr] & & 1 \\
& & & \prod_{v'}H^1(K'_{v'},\overline{G}) \ar'[r][rr] & & \prod_{v'}H^2(K'_{v'},Z)\\
\prod_{v'}H^1(K',_{\zeta'}\widetilde{G})\ar[rr]_{\gamma'_{V',\zeta'}}& &\prod_{v'}H^1(K',_{\zeta'}\overline{G}) \ar[rr]_{\delta'_{V',\zeta'}}\ar[ur]_-{\tau_{V',\zeta'}}& & \prod_{v'}H^2(K',Z) \ar[ur]_-{\tau_{V',\eta'}} }
\end{equation}
As usual $\tau_{V',\zeta'}=\prod_{v'}\tau_{v',\zeta'}$ and $\tau_{V',\eta'}=\prod_{v'}\tau_{v',\eta'}$ where $\tau_{v,\zeta'},\tau_{v',\eta'}$ are the local counterparts to $\tau_{\zeta'}$ and $\tau_{\eta'}$. The two unlabeled arrows in the backmost face of the cube are the maps $\delta_{V'}$ and $\iota_{\overline{G},K'}$. The two unlabeled arrows in the front face of the cube are the maps $\iota_{Z,K',\zeta'}$ and $\delta'_{\zeta'}$. The existence of $\Phi_{\widetilde{\Ad}}$ follows from Lemma \ref{LocalEquivCoho} and the bijectivity of $\iota_{\widetilde{G},K},\iota_{\widetilde{G},K'}$ (see \cite[Thm 6.6]{PR}). Explicitly we have $\Phi_{\widetilde{\Ad}} = \iota_{\widetilde{G},K'}^{-1} \circ \vp_0 \circ \iota_{\widetilde{G},K}$.

Returning to the Galois cohomological diagram (\ref{E:TheCube}), the bijectivity of $\tau_{V',\zeta'}$ implies that there exists $(\theta'_{v'}) \in \prod_{v'}H^1(K',_{\zeta'}\overline{G})$ such that $\tau_{V',\zeta'}((\theta'_{v'}))=(\xi'_{v'})$. Via the commutativity of the bottom face of the cube, $\delta'_{V',\zeta'}(\theta'_{v'})=(1_{v'})$ where $(1_{v'})$ denotes the trivial cocycle in $\prod_{v'} H^2(K,Z)$. Consequently there exists $(\mu'_{v'})\in\prod_{v'}H^1(K_{v'},_{\zeta'}\widetilde{G})$ such that $\gamma'_{V',\zeta'}((\mu'_{v'})) = (\theta'_{v'})$. As $\iota_{\widetilde{G},K',\zeta'}$ is bijective, there exists $\mu'\in H^1(K',_{\zeta'}\widetilde{G})$ such that $\iota_{\widetilde{G},K',\zeta'}(\mu')=(\mu'_{v'})$. Setting $\xi' = \tau_{\zeta'}(\gamma'_{\zeta'}(\mu'))$, we obtain that $\iota_{\overline{G},K'}(\xi')=(\xi'_{v'})$ and $\delta'(\xi') = \eta'$.
This is the desired cocycle.
\end{proof}
\end{proof}

\subsection{Proof of Theorem \ref{GC-MaximalBiject}}

As $H^1(K,\overline{G})$ bijects with inner forms of $G$, to each $[\xi]$ we may associate (up to equivalence) an inner twist of $G$ which we denote $_\xi G$. After identification, $\Phi_{\Ad}$ takes the $K$-isomorphism class of $_\xi G$ to the $K'$--isomorphism class of $_{\xi'}G'$. We now show that there is a bijection between maximal arithmetic lattices of a fixed group $H$ in the class of $_\xi G$ and a fixed group $H'$ in the class of $_{\xi'}G'$. This is done using the characterization of maximal arithmetic lattices furnished by Borel--Prasad \cite[Prop 1.4]{BP}.

\begin{thm}[Borel--Prasad, \cite{BP}]\label{BPThm}
Let $H$ be an absolutely almost simple $K$-group and $\Gamma$ a maximal arithmetic subgroup of $\prod_{v\in V^K_\infty}H(K_v)$ (where we assume this product has no compact factors). Further, let $\Lambda$ be the inverse image of $\Gamma\cap i(\widetilde{H}(K))$ in $\widetilde{H}(K)$ where $i:\widetilde{H}\to H$ is a central isogeny from the simply connected form $\widetilde{H}$. Then
\begin{itemize}
\item For $v\notin V_\infty^K$, the closure of $\Lambda$ in $\widetilde{H}(K_v)$ is a parahoric subgroup of $\widetilde{H}(K_v)$.
\item $\Lambda=\widetilde{H}(K)\cap\prod_{v\in V^K_f} P_v$.
\item $\Gamma$ is the normalizer of $i(\Lambda)$ in $\prod_{v\in V^K_\infty}H(K_v)$.
\end{itemize}
\end{thm}

Recall our earlier notation that $v'=\Phi_{\V}(v)$ and we use the notation $H_v$ to denote $H(K_v)$, with similar notation for $H'$. Theorem \ref{BPThm} implies that to construct a bijection between maximal arithmetic lattices, it suffices to construct a bijection between principal arithmetic lattices, as taking the image under $i$ and then normalizers yields the result. $K_v$ and $K'_{v'}$ are isomorphic via $\Phi_{\V}$ and so $\widetilde{H}_v$ and $\widetilde{H}'_{v'}$ are isomorphic. Consequently, there is an isomorphism between their corresponding buildings $\mathcal{B}_{K_v}(\widetilde{H}_v)$ and $\mathcal{B}_{K'_{v'}}(\widetilde{H}'_{v'})$ denoted $\Phi_{\Build}$, which preserves the chamber structure. As parahoric subgroups arise as the stabilizers of facets in the building, this isomorphism induces a bijection between parahoric subgroups of $\widetilde{H}_v$ and $\widetilde{H}'_{v'}$.
Indeed, any facet $F\in\mathcal{B}_{K_v}(\widetilde{H}_v)$ is mapped to a corresponding facet $\Phi_{\Build}(F)\in\mathcal{B}_{K'_{v'}}(\widetilde{H}'_{v})$ and vice versa.
Consequently we have an induced bijection, $\Phi_{\Para,v}$, given by $\Phi_{\Para,v}(\widetilde{H}_v^F)=(\widetilde{H}'_{v'})^{\Phi_{\Build}(F)}$. Via $\Phi_{\Para,v}$ we have an induced bijection $\Phi_{\Para}$ between coherent collections of parahorics taking $\textbf{P}=\prod_{v\in V^K_f}P_v$ to a corresponding coherent collection $\textbf{P}'=\prod_{v'\in V^{K'}_f}\Phi_{\Para,v}(P_v)$. By Theorem \ref{BPThm}, $\Phi_{\Para}$ gives rise to a bijection between principal arithmetic lattices associated to $\Lambda$ and $\Lambda'$ by taking intersections with the $K$ and $K'$ points of $\widetilde{H}$ and $\widetilde{H}'$, respectively. By taking normalizers we get the desired result, Theorem \ref{GC-MaximalBiject}.



 \end{document}